\documentclass[12pt]{article}
\usepackage{amsthm, amsmath,amssymb,amsfonts}
\usepackage{graphicx}
\newtheorem{theorem}{Theorem}[section]
\newtheorem{proposition}[theorem]{Proposition}
\newtheorem{corollary}[theorem]{Corollary}
\newtheorem{lemma}[theorem]{Lemma}
\newtheorem{remark}[theorem]{Remark}
\newtheorem{definition}[theorem]{Definition}
\newtheorem{example}[theorem]{Example}
\begin{document}
\title{
Generalizing circles over algebraic extensions}
\author{T. Recio\footnote{The authors are partially supported by the project
MTM2005-08690-CO2-01/02 ``Ministerio de Educaci\'on y Ciencia".}
\and J.R. Sendra\footnotemark[1] \footnote{Partially
supported by CAM-UAH2005/053 ``Direcci\'on General de Universidades
de la Consejer\'{\i}a de Educaci\'on de la CAM y la Universidad de
Alcal\'a".} \and L.F. Tabera\footnotemark[1] \footnote{L.F.Tabera
also supported by a FPU research grant.} \and C.
Villarino\footnotemark[1] \footnotemark[2]}
\maketitle

\begin{abstract}
This paper deals with a family of spatial rational curves that were introduced
in
 \cite{ARS-2}, under the name of {\it hypercircles}, as an algorithmic
cornerstone tool  in the context of improving the rational
parametrization (simplifying the coefficients of the rational
functions, when possible) of algebraic varieties. A real circle can
be defined as the image of the real axis under a Moebius
transformation in the complex field. Likewise, and roughly speaking,
a hypercircle can be defined as the image of a line (``the
${\mathbb{K}}$-axis") in a $n$-degree finite algebraic extension
$\mathbb{K}(\alpha)\thickapprox\mathbb{K}^n$ under the
transformation
$\frac{at+b}{ct+d}:\mathbb{K}(\alpha)\rightarrow\mathbb{K}(\alpha)$.

The aim of this article is to extend, to the case of hypercircles,
some of the specific properties of circles. We show that
hypercircles are precisely, via $\mathbb{K}$-projective
transformations, the rational normal curve of a suitable degree. We
also obtain a complete description of the points at infinity of
these curves (generalizing the cyclic structure at infinity of
circles). We characterize hypercircles as those curves of degree
equal to the dimension of the ambient affine space and with
infinitely many ${\mathbb{K}}$-rational points, passing through
these points at infinity. Moreover, we give explicit formulae for
the parametrization and implicitation of hypercircles. Besides the
intrinsic interest of this very special family of curves, the
understanding of its properties has a direct application to the
simplification of parametrizations problem, as shown in the last
section.
\end{abstract}

\section{Introduction}
The problem of obtaining a real parametrization of a rational planar
curve given by a complex parametrization has been studied --from an
algorithmic point of view-- in \cite{RS-Real}. There, the problem is
reduced to determining that a certain curve obtained after
manipulating the given parametrization is a real  line or a real
circle. From a real parametrization of this circle (or line), a real
parametrization of the original curve is then achieved. This
auxiliary circle is found by an analogous to Weil descente's method
\cite{Weil} applied to the complex parametrization of the originally
given curve. In \cite{ARS-2}, the same approach has been extended to
the general case of planar or spatial rational curves $\mathcal{C}$
given by a parametrization over $\mathbb{K}(\alpha)$, where $\alpha$
is an algebraic element over $\mathbb{K}$. In order to obtain,
whenever possible, a parametrization over $\mathbb{K}$ of
$\mathcal{C}$, another rational curve, with remarkable properties,
is associated to $\mathcal{C}$. In \cite{ARS-2} it is shown that
this associated curve is, in the relevant cases, a generalization of
a circle, in the sense we will discuss below, deserving to be named
hypercircle.

The simplest hypercircles should be the circles themselves. We can think of the
real plane as the field of  complex numbers $\mathbb{C}$, an algebraic
extension of the reals $\mathbb{R}$ of degree 2. Analogously, we can
consider a characteristic zero base field $\mathbb{K}$ and an algebraic
extension of degree $n$, $\mathbb{K}(\alpha)$. Let us identify
$\mathbb{K}(\alpha)$ as the vector space $\mathbb{K}^n$, via the choice of
a suitable base, such as the one given by the powers of $\alpha$. This is the
framework in which hypercircles are defined.

Now let us look to the different, equivalent,  ways of defining a common circle
on the real plane, with the purpose of taking the most convenient one for
generalization. The first definition of a circle is the set of points in the
real
plane that are equidistant from a fixed point. This approach does not extend
well to more general algebraic extensions, because we do not have an
immediate notion of metric over $\mathbb{K}^n$. On the other hand,
algebraically, a real planar circle is a conic such that its homogeneous degree
two form is $x^2+y^2$ and such that it contains an infinite number of real
points. Even if we will prove in Section 6 that we can show an analogous
definition for a hypercircle, this is not an operative way to start defining
them.

Finally, from another point of view, we see that circles are real
rational curves. This means that there are two real rational
functions $(\phi_1(t),\phi_2(t))$ whose image cover almost all the
points of the circle. For instance, the circle $x^2+y^2=1$ is
parametrized by $\phi(t)=(\frac{t^2-1}{t^2+1},\frac{2t}{t^2+1})$.
Every proper (almost one-to-one \cite{SW-Symbolic}) rational
parametrization of a circle verifies that
$\phi_1(t)+i\phi_2(t)=\frac{at+b}{ct+d}\in\mathbb{C}(t)\setminus
\mathbb{C}$, which defines a conformal mapping
$u:\mathbb{C}\rightarrow\mathbb{C}$. Moreover, if we identify
$\mathbb{C}$ with $\mathbb{R}^2$, the image of the real axis $(t,0)$
under $u$ is exactly the circle parametrized by $\phi(t)$.
Conversely, let $u(t)=\frac{at+b}{ct+d}\in\mathbb{C}(t)$ be a unit
of the near-ring $\mathbb{C}(t)$ under the composition operator (see
\cite{Near}). If $c\neq 0$ and $d/c\notin \mathbb{R}$ then, the
closure of the image by $u$ of the real axis is a circle. Otherwise,
it is a line. This method to construct circles generalizes easily to
algebraic extensions. Namely, let $u(t)=\frac{at+b}{ct+d}$ be a unit
of $\mathbb{K}(\alpha)(t)$ (i.e. verifying that $ad-bc \neq 0$). Let
us identify $\mathbb{K}(\alpha)$ with $\mathbb{K}^n$ and let $u$ be
the map
\[\begin{matrix}
u:&\mathbb{K}(\alpha)\approx \mathbb{K}^{n}&\rightarrow&
\mathbb{K}(\alpha)\approx {\mathbb{K}}^{n}\\
&t&\mapsto& u(t)\end{matrix}.\]
Then, the Zariski-closure of the image of the axis $(t,0,\ldots,0)$ under the
map $u$ is a rational curve in $\mathbb{K}^n$. These curves are,
by definition, our hypercircles.

Roughly speaking,  it happens (see \cite{ARS-2}) that a parametrization over
$\mathbb{K}$ of  the hypercircle associated to a given rational curve
$\mathcal{C}$ (whose parametrization we want to simplify) can be used to get
--in a straightforward manner-- a parametrization of $\mathcal{C}$  over
$\mathbb{K}$. As pointed in \cite{ARS-2}, it seems that, due to the geometric
properties of hypercircles, it is algorithmically simpler to obtain such
parametrization for this type of curves than it is for $\mathcal C$. In fact, it
is
shown in \cite{issac04} how to get this in some cases. Therefore, the
reparametrization problem is behind our increasing interest in the study of
hypercircles on its own.

The structure of this paper is as follows. In Section 2 we formally
introduce the notion of hypercircle. We study the influence on a
hypercircle when adding and multiplying  the  defining unit $u(t)$
by elements of $\mathbb{K}(\alpha)$, reducing the affine
classification of hypercircles to those defined by some simpler
units. Next we characterize the units associated to lines. In
Section 3 we show how to transform, projectively, a
hypercircle into the {\it rational normal curve}
(see \cite{harris}). From this, we derive the main geometric
properties of hypercircles (smoothness, degree, affine equivalence,
etc.) and we reduce the study of hypercircles to the subclass of
{\it primitive hypercircles} (See Definition
\ref{def_hipercirculo_primitivo}).  In Section 4 the behavior of
hypercircles at infinity is analyzed, showing its precise and rich
structure. In Section 5, exploiting the stated geometric features,
we present {\it ad hoc} parametrization and implicitization methods
for hypercircles. In Section 6 we characterize hypercircles among
curves of degree equal to the dimension of the ambient affine space,
passing through the prescribed points of infinity described in
Section 4 and having infinitely many rational points.
Finally, Section 7 is devoted to show how the insight gained
throughout  this paper can be applied to derive heuristics for
solving the problem of simplifying the parametrization of curves
with coefficients involving algebraic elements.

Throughout this paper the following notation and terminology will be used.

\begin{itemize}
\item $\mathbb{K}$ will be a field of characteristic zero,
$\mathbb{K}\subseteq \mathbb{L}$ a finite algebraic extension of degree
$n$ and $\mathbb{F}$ the algebraic closure of $\mathbb{K}$.
\item $\alpha$ will be a primitive element of $\mathbb{L}$ over $\mathbb{K}$.
\item $u(t)$ will be a unit under composition of $\mathbb{L}(t)$. That is,
$u(t)=\frac{at+b}{ct+d}$ with $ad-bc\neq 0$.
Its inverse $\frac{-dt+b}{ct-a}$ is denoted by $v(t)$.
\item For $u(t)=\frac{at+b}{ct+d}$ and $c\neq 0$,
$M(t)=t^r+k_{r-1}t^{r-1}+\cdots+k_0\in\mathbb{K}[t]$ denotes the
minimal polynomial of $-d/c$ over $\mathbb{K}$.
\item We will denote as $m(t)$ the polynomial obtained by dividing $M(t)$ by
$ct+d$. That is, $\displaystyle{m(t)=\frac{M(t)}{ct+d}=
l_{r-1}t^{r-1}+l_{r-2}t^{r-2}+\cdots+l_0\in\mathbb{L}[t]}$.
\item Sometimes we will represent $u(t)$ as
\[\displaystyle{u(t)=\frac{(at+b)m(t)}{M(t)}=\frac{p_0(t)+p_1(t)\alpha+\cdots+p_
{n-1}(t)\alpha^{n-1}}{M(t)}},\]
where $p_i(t)\in \mathbb{K}[t]$.
\item By $\{\sigma_1=Id,\sigma_2,\ldots,\sigma_s\}$, $s\geq n$ we will denote
the group of $\mathbb{K}$-automorphisms of the normal closure of
$\mathbb{K}\subseteq \mathbb{L}$.
\item We will represent by $\{\alpha_1=\alpha,\ldots,\alpha_n\}$ the
conjugates of $\alpha$. We assume without loss of generality that
$\sigma_i(\alpha)=\alpha_i$ for $i=1,\ldots,n$.
\end{itemize}

\section{Definition and First Properties}

In this section we begin with the formal definition of a hypercircle.
\begin{definition}\label{def_hipercirculo}
Let $u(t)$ be a unit in $\mathbb{L}(t)$, where
${\mathbb{L}}=\mathbb{K}(\alpha)$. Let
\[u(t)=\sum_{i=0}^{n-1}\phi_i(t)\alpha^i\]
where $\phi_i(t)\in\mathbb{K}(t)$, for $i=0,\ldots,n-1$. The
$\alpha$-hypercircle $\mathcal{U}$ generated by $u(t)$ is the rational curve
in $\mathbb{F}^n$ parametrized by $\phi(t)=(\phi_0(t),\ldots,\phi_{n-1}(t)).$
\end{definition}

Observe that the expansion of $u(t)$ in powers of $\alpha$ is
unique, because $\{1,\alpha,\ldots,\alpha^{n-1}\}$ is a basis of
$\mathbb{K}(\alpha)(t)$ as a $\mathbb{K}(t)-$vector space. The
parametrization can be obtained by rationalizing the denominator as
follows: suppose given the unit $u(t)=\frac{at+b}{ct+d}$, $c\neq 0$
(remark that, if $c=0$, it is straightforward to obtain $\phi(t)$),
and the extension $\mathbb{K}\subseteq \mathbb{K}(\alpha)$. Let
$M(t)$ be the minimal polynomial of $-d/c$ over $\mathbb{K}$.
Compute the quotient $m(t)=\frac{M(t)}{ct+d}\in
\mathbb{K}(\alpha)[t]$ and develop the unit as
\[\frac{{a}t+b}{ct+d}=\frac{(at+b)m(t)}{M(t)}=\frac{p_0(t)+p_1(t)\alpha+
\cdots+p_{n-1}(t)\alpha^{n-1}}{M(t)}\] where $p_i(t)\in
\mathbb{K}[t]$. From this,
$\phi(t)=\left(\frac{p_0(t)}{M(t)},\ldots,
\frac{p_{n-1}(t)}{M(t)}\right)$ is the parametrization associated to
$u(t)$. Remark that $\gcd(p_0(t),\ldots, p_{n-1}(t), M(t))=1$.
Moreover, it is clear that
$\mathbb{F}(\phi_0(t),\ldots,\phi_{n-1}(t))=\mathbb{F}(t)$. So this
parametrization is proper in $\mathbb{F}$, and it follows from the
results in \cite{Gutierrez-Recio-near-separated} that also
$\mathbb{K}(\phi_0(t),\ldots,\phi_{n-1}(t))=\mathbb{K}(t)$.

\begin{figure}\begin{center}
\includegraphics[width=6cm]{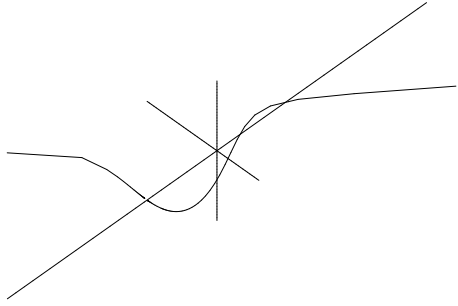}\caption{A hypercircle in
$\mathbb{R}^3$}\label{fig:Figura_grado3_una_rama}
\end{center}\end{figure}

\begin{example}\label{ejemplo2}
Let us consider the algebraic extension
$\mathbb{Q}\subseteq\mathbb{Q}(\alpha)$, where
$\alpha^3+2\alpha+2=0$. The unit $\frac{t-\alpha}{t+\alpha}$ has an
associated hypercircle parametriced by
\[\phi(t)=\left(\frac{t^3+2t+2}{t^3+2t-2},\frac{-2t^2}{t^3+2t-2},
\frac{2t}{t^3+2t-2}\right)\]
A picture of the spatial real curve is shown in Figure
\ref{fig:Figura_grado3_una_rama}
\end{example}

As it stands, the definition of a hypercircle $\mathcal{U}$ depends on a given
unit $u(t) \in \mathbb{L}(t)$  and on a primitive generator $\alpha$ of an
algebraic extension $\mathbb{L}$.  In what follows we will analyze the effect
on $\mathcal{U}$ when varying some of these items, searching for a simple
representation of a hypercircle to ease studying its geometry.

First notice that, given a unit $u(t) \in \mathbb{L}(t)$ and two different
primitive elements $\alpha$ and $\beta$ of the extension
$\mathbb{K}\subseteq\mathbb{L}$, we can expand the unit in two different
ways $u(t)=\sum_{i=0}^{n-1}\alpha^i\phi_i(t)=
\sum_{i=0}^{n-1}\beta^i\psi_i(t)$. The hypercircles
$\mathcal{U}_\alpha\simeq(\phi_0(t),\ldots,\phi_{n-1}(t))$ and
$\mathcal{U}_\beta\simeq(\psi_0(t),\ldots,\psi_{n-1}(t))$ generated by $u(t)$
are different curves in $\mathbb{F}^n$, see Example \ref{ejemplo1}.
Nevertheless, let $\mathcal{A}\in\mathcal{M}_{n\times n}(\mathbb{K})$ be
the matrix of change of basis from $\{1,\alpha,\ldots,\alpha^{n-1}\}$ to
$\{1,\beta,\ldots,\beta^{n-1}\}$. Then,
$\mathcal{A}(\phi_0(t),\ldots,\phi_{n-1}(t))^t=
(\psi_0(t),\ldots,\psi_{n-1}(t))^t$. That is, it carries one of the curve onto
the other.
Thus, $\mathcal{U}_\alpha$ and $\mathcal{U}_{\beta}$ are related by the
affine transformation  induced by the change of basis and, so, they share many
important geometric properties.

In the sequel, if there is no confusion about the algebraic extension and the
primitive element, we will simply call $\mathcal{U}$ a hypercircle.

\begin{example}\label{ejemplo1}
Let us consider the algebraic extension
$\mathbb{Q}\subseteq\mathbb{Q}(\alpha)$, where $\alpha^4+1=0$. Let us
take the unit $u(t)=\frac{t-\alpha}{t+\alpha}$. By normalizing $u(t)$, we
obtain the parametrization $\phi(t)$ associated to $u(t)$:
\[\phi(t)=\left(\frac{t^4-1}{t^4+1},
\frac{-2t^3}{t^4+1},\frac{2t^2}{t^4+1}, \frac{-2t}{t^4+1}\right).\]
This hypercircle ${\mathcal{U}}_\alpha$ is the zero set of
$\{X_1X_2-X_3X_0-X_3,X_1^2+X_3^2-2X_2,X_1X_0+X_2X_3-X_1,X_0^2+X_3X_1-1\}$.
Now, we take $\beta=\alpha^3+1$, instead of $\alpha$, as the
primitive element of $\mathbb{Q}(\alpha)=\mathbb{Q}(\beta)$. The
same unit $u(t)$ generates the $\beta$-hypercircle
${\mathcal{U}}_\beta$ parametrized by
\[\psi(t)=\left(\frac{t^4+2t^3-2t^2+2t-1}{t^4+1},\frac{-6t^3+4t^2-2t}{t^4+1},
\frac{6t^3-2t^2}{t^4+1},\frac{-2t^3}{t^4+1}\right),\]
which is different to ${\mathcal{U}}_\alpha$; note that $\psi(1)=(1,-2,2,-1)$
that does not satisfy the  equation $X_0^2+X_3X_1-1=0$ of
${\mathcal{U}}_\alpha$.
\end{example}

On the other hand it is well known that a given parametric curve can
be parametrized over a given field $\mathbb{S}$ by different proper
parametrizations, precisely, those obtained by composing to the
right a given proper parametrization by a unit in $\mathbb{S}(t)$.
In this way, we have a bijection between $\alpha$-hypercircles and
the equivalence classes of units of $\mathbb{K}(\alpha)(t)$ under
the equivalence relation {``\it $u\sim v$ iff $u(t)=v(\tau(t))$ for
a unit $\tau(t)\in\mathbb{K}(t)$''} (fixing the correspondence,
between a unit in $\mathbb{K}(\alpha)(t)$ and a hypercircle,
by means of the expansion of the unit in terms of powers of $\alpha$).

More interesting is to analyze, on a hypercircle defined by a unit $u(t)$, the
effect of composing it to the left with another unit
$\tau(t) \in \mathbb{K}(\alpha)(t)$, that is, of getting $\tau(u(t))$.
For instance, $\tau(t)$ could be $\tau(t)=t+\lambda$ or $\tau(t)=\lambda t$,
or $\tau(t)=1/ t$, with $\lambda \in \mathbb{K}(\alpha)^{\ast}$. Every unit is
a sequence of compositions of these three simpler cases,  for instance, when
$c\neq 0$, we have
\[t\longmapsto ct \longmapsto ct+d \longmapsto \frac{1}{ct+d}\longmapsto
\frac{bc-ad}{c}\frac{1}{ct+d}  \longmapsto \]
\[\longmapsto \frac{a}{c}+
\frac{bc-ad}{c}\frac{1}{ct+d}=\frac{at+b}{ct+d}=u(t).\]
Therefore, studying their independent effect is all we need to understand
completely the behavior of a hypercircle under left composition by units.

For circles, adding a complex number to the unit that defines the circle
correspond to a translation of the circle. Multiplying it by a complex number
acts as the composition of a rotation and a dilation. And the application
$\tau(t)=1/ t$ gives an inversion. The following lemma analyzes what happens
in the general case.

\begin{lemma}\label{lema-traslacion}
Let $\mathcal{U}$ be the $\alpha$-hypercircle generated by $u(t)$,
and
$\displaystyle{\lambda=\sum_{i=0}^{n-1}\lambda_i\alpha^i\in\mathbb{K}(\alpha)^{*
}}$,
where $\lambda_i\in\mathbb{K}$. Then,
\begin{enumerate}
\item $\lambda+u(t)$ is a unit generating the hypercircle obtained from
$\mathcal{U}$ by the translation of vector
$(\lambda_0,\ldots,\lambda_{n-1})$.
\item $\lambda u(t)$ is a unit generating the hypercircle obtained from
$\mathcal{U}$ by the affine transformation over $\mathbb{K}$ given by the
matrix of change of basis from
$\mathcal{B}^{\star}=\{\lambda,\lambda\alpha,\ldots,\lambda\alpha^{n-1} \}$ to
$\mathcal{B}=\{1,\alpha,\ldots,\alpha^{n-1}\}$.
\end{enumerate}
\end{lemma}

\begin{proof} To prove (1), let
$\phi(t)=(\phi_0(t),\ldots,\phi_{n-1}(t))\in\mathbb{K}(t)^n$ be the
parametrization of $\mathcal{U}$ obtained from $u(t)$.
Then, $\lambda+u(t)=\sum_{i=0}^{n-1}(\lambda_i+\phi_i(t))\alpha^i$
generates the hypercircle parametrized by
$(\lambda_0+\phi_0(t),\ldots,\lambda_{n-1}+\phi_{n-1}(t))\in\mathbb{K}(t)^n$,
which is the translation of $\mathcal{U}$ of vector
$(\lambda_0,\ldots,\lambda_{n-1})$. For the second assertion, let
$\phi^{\star}(t)\in\mathbb{K}(t)^ n$ be the parametrization of the
hypercircle associated to the unit $\lambda u(t)$. The rational coordinates
$\phi^{\star}_i(t)$ of $\phi^{\star}(t)$ are obtained from the matrix
$\mathcal{A}=(a_{i,j})\in\mathcal{M}_{n\times n}(\mathbb{K})$ of change of basis
from $\mathcal{B}^{\star}$ to $\mathcal{B}$, for $i,j=0,\ldots,n-1$. Indeed,
\[\lambda u(t)=\sum_{i=0}^{n-1}\phi_i(t)\lambda\alpha^i=
\sum_{i=0}^{n-1}\phi_i(t)\left(\sum_{j=0}^{n-1}a_{ji}\alpha^{j}\right) =
\sum_{j=0}^{n-1}\left(\sum_{i=0}^{n-1}a_{ji}\phi_i(t)\right)\alpha^{j}.\]
Then $\phi^{\star}(t)^{t}=\mathcal{A}\,\phi(t)^{t}$.
\end{proof}

Finally, the following lemma uses the previous results to transform affinely one
hypercircle into another one whose unit is simpler.

\begin{lemma}\label{lema-forma-reducida}
Let $u(t)=\frac{at+b}{ct+d}$ be a unit and $\mathcal{U}$ its associated
hypercircle.
\begin{enumerate}
\item If $c=0$ then $\mathcal{U}$ is affinely equivalent over $\mathbb{K}$
to the line generated by $u^\star(t)=t$.
\item If $c\neq 0$ then $\mathcal{U}$ is affinely equivalent over
$\mathbb{K}$ to the hypercircle $\mathcal{U}^\star$ generated by
$u^{\star}(t)=\frac{1}{t+d/c}$
\end{enumerate}
\end{lemma}
\begin{proof}
This lemma follows from Lemma \ref{lema-traslacion}, taking into account that
$u(t)$ is obtained from $u^{\star}(t)$ by the following composition:
\[u^{\star}(t)\mapsto \lambda_1 u^{\star}(t) \mapsto
\lambda_1u^{\star}(t)+\lambda_2=u(t)\]
with suitable $\lambda_1, \lambda_2, u^{\star}$. If $c=0$, then
$\lambda_1=\frac{a}{d}\neq 0$ and $\lambda_2=\frac{b}{d}$ for
$u^{\star}(t)=t$. Analogously, if $c\neq 0$, then $u(t)$ is obtained from
$u^{\star}(t)=\frac{1}{t+d/c}$ taking $\lambda_1=\frac{bc-ad}{c^2}\neq 0$
and $\lambda_2=\frac{a}{c}$.
\end{proof}

Therefore the (affine) geometry of hypercircles can be reduced to those
generated by a unit of type $\frac{1}{t+d}$ (then we say the unit is in
{\it reduced form}). The simplest hypercircle of this kind is given by
$\frac{1}{t+d}$, when $d\in \mathbb{K}$. It is the line parametrized by
$(\frac{1}{t+d},0,\ldots,0)$. In the complex case, Moebius transformations
defining lines are precisely those given either by a polynomial unit in $t$
(i.e. a
unit without $t$ at the denominator) or by a unit such that the root of the
denominator is in $\mathbb{R}$. The same property holds for hypercircles.

\begin{theorem}\label{teo-rectas}
Let $\mathcal{U}$ be the $\alpha$-hypercircle associated to $u(t)$. Then, the
following statements are equivalent:
\begin{enumerate}
\item $\mathcal{U}$ is a line.
\item $\mathcal{U}$ is associated to a polynomial unit.
\item The root of the denominator of every non polynomial unit generating
$\mathcal{U}$ belongs to $\mathbb{K}$.
\item $\mathcal{U}$ is polynomially parametrizable (over $\mathbb{F}$).
\item $\mathcal{U}$ has one and only one branch (over $\mathbb{F}$ ) at
infinity.
\item $\mathcal{U}$ is polynomially parametrizable over $\mathbb{K}$.
\item $\mathcal{U}$ has one and only one branch (over $\mathbb{K}$ ) at
infinity.
\end{enumerate}
\end{theorem}

\begin{proof}
(1) $\Leftrightarrow$ (2). By definition, we know that hypercircles have a
parametrization over $\mathbb{K}$. Thus, if $\mathcal{U}$ is a line, it can be
parametrized as $(a_0t+b_0,\ldots,a_{n-1}t+b_{n-1})$, where
$a_i\,,\,b_i\,\in\,{\mathbb{K}}$. Therefore,
$u(t)=\left(\sum_{i=0}^{n-1}\,a_i\alpha^i\right)\,t+\sum_{i=0}^{n-1}\,
b_i\alpha^i$
is a polynomial unit associated to $\mathcal{U}$. Conversely, let $u(t)=at+b \in
\mathbb{L}(t)$,
$a\neq 0$, be a polynomial unit associated to $\mathcal{U}$.  Then $\mathcal{U}$
is the
line parametrized by
$\mathcal{P}(t)=(a_0t+b_0,\ldots,a_{n-1}t+b_{n-1})\in\mathbb{K}[t]^n,$
where $a=\sum_{i=0}^{n-1}\,a_i\,\alpha^i$ and
$b=\sum_{i=0}^{n-1}\,b_i\,\alpha^i$.

(2) $\Leftrightarrow$ (3). Let $u(t)=at+b$ be a polynomial unit associated to
$\mathcal{U}$, and let $u^{\star}(t)$ be another non polynomial unit associated
to
$\mathcal{U}$. Then, $u^{\star}(t)=u(\tau(t))$, where $\tau(t)$ is a unit of
$\mathbb{K}\,(t)$. Therefore, the root of $u^{\star}(t)$ belongs to
$\mathbb{K}$. Conversely, by Lemma \ref{lema-forma-reducida}, (3) implies
(1), and we know that (1) implies (2).

(3) $\Leftrightarrow$ (4). Indeed, (3) implies (2) and therefore
(4). Conversely, let $u(t)$ be a non-polynomial unit generating
$\mathcal{U}$, and let $\phi(t)=(\phi_i)_{i=1,\ldots,n}
\in\mathbb{K}(t)^n$ be the associated parametrization of
$\mathcal{U}$. Then, $\phi(t)$ is proper,
$\phi_i(t)=\frac{p_i(t)}{M(t)}$ with $\deg(p_i)\leq \deg(M)$ and
$\gcd(p_0(t)\dots p_{n-1}(t), M(t))=1$. Thus, the fact that
$\mathcal{U}$ admits a polynomial parametrization, implies, by
Abhyankar-Manocha-Canny's criterion of polynomiality (see
\cite{manocha}), that the denominator $M(t)$ is either constant or
has only one root. Now, $M(t)$ can not be constant, since it is a
minimal polynomial. Thus, $M$ has only one root, and since it is
irreducible, it must be linear. Moreover, since $M\in
{\mathbb{K}}[t]$, its root is an element in $\mathbb{K}$.

(4) $\Leftrightarrow$ (5) This is, again, the geometric version of
Abhyankar-Manocha-Canny's criterion. Same for (6) $\Leftrightarrow$ (7).

(4) $\Leftrightarrow$ (6) Obviously (6) implies (4). Conversely,  if we have a
polynomial parametrization over $\mathbb{F}$, it happens \cite{AR} that any
proper parametrization must be either polynomial or in all its components the
degree of the numerator must be smaller or equal than the degree of the
denominator and, then, this denominator has only one single root over
$\mathbb{F}$. So, since the parametrization $\phi(t)$ induced by the unit is
proper, and by hypothesis  $\mathcal{U}$ is polynomial, then $\phi(t)$ must be
either polynomial (in which case we are done because $\phi(t)$ is over
$\mathbb{K}$) or its denominator $M(t)$ has a single root $a\in \mathbb{F}$.
Now, reasoning as above one gets that $a\in\mathbb{K}$. So, a change of
parameter, such as $t \mapsto \frac {1+as}{s}$ turns $\phi(t)$ into a
$\mathbb{K}$-polynomial parametrization.
\end{proof}

As a corollary of this theorem, we observe that a parabola  can never be  a
hypercircle, since it is polynomially parametrizable, but it is not a line.
Nevertheless, it is easy to check that the other irreducible conics are indeed
hypercircles for certain algebraic extensions of degree 2.

\section{Main Geometric Properties.}
This section is devoted to the analysis on the main geometric
properties of hypercircles. The key idea, when not dealing with
lines, will be to use the reduction to units of the form
$u(t)=\frac{1}{t+d}$, where $d\notin\mathbb{K}$ (see Lemma
\ref{lema-forma-reducida}).

\begin{theorem}\label{teo-estructura-hipercirculo}
Let $\mathcal{U}$ be the $\alpha$-hypercircle associated to the unit
$u(t)=\frac{at+b}{t+d}\in\mathbb{K}(\alpha)(t)$ and let
$r=[\mathbb{K}(-d):\mathbb{K}]$. Then,
\begin{enumerate}
\item there exists an affine transformation
$\chi:\mathbb{F}^{n}\longrightarrow \mathbb{F}^{n}$ defined over
$\mathbb{K}$ such that the curve $\chi(\mathcal{U})$ is parametrized by
\[\widetilde{\chi}(t)=\left(\frac{1}{M(t)},\frac{t}{M(t)},\ldots,\frac{t^{r-1}}{
M(t)},0,\ldots,0\right).\]
\item there exists a projective transformation
$ \rho:\mathbb{P(F)}^{n}\longrightarrow \mathbb{P(F)}^{n}$, defined
over $\mathbb{K}$, such that the curve $\rho(\mathcal{U})$ is the
rational normal curve of degree $r$ in $\mathbb{P(F)}^n$,
parametrized by
\[\widetilde{\rho}(t:s)=[s^r:s^{r-1}t:\cdots:st^{r-1}:t^r:0:\cdots:0].\]
\end{enumerate}
\end{theorem}
\begin{proof} For the case of lines the result is trivial. By Lemma
\ref{lema-forma-reducida}, we can consider that $\mathcal{U}$ is the
hypercircle associated to $u(t)=\frac{1}{t+d}$ and $r\geq 2$. Let
$M(t)=t^r+k_{r-1}t^{k-1}+\cdots+k_0\in
\mathbb{K}[t],m(t)=\sum_{i=0}^{r-1} l_{i}t^i\in {\mathbb{L}}[t]$, as
indicated in Section 1 and, since the numerator of $u(t)$ is 1, it
holds that $m(t)=\sum_{i=0}^{n-1} p_i(t)\alpha^i$, $p_i(t)\in
\mathbb{K}[t]$. Also, note that both $M(t)$ and the denominator of
$u(t)$ are monic, and hence $l_{r-1}=1$. First of all, we prove that
there are  exactly $r$ polynomials in $\{p_i(t),\
i=0,\ldots,n-1\}\subset \mathbb{K}[t]$ being linearly independent.
For this purpose, we observe that the coefficients of $m(t)$,
$\{1,l_{r-2},\ldots,l_0\}\subset \mathbb{L}$, are linearly
independent over $\mathbb{K}$. Indeed, from the equality
$M(t)=(t+d)m(t)$, one has that
$l_{r-i}=(-d)^{i-1}+(-d)^{i-2}k_{r-1}+\cdots+k_{r-i+1}$, for
$i=2,\ldots,r.$ So, $\{1,l_{r-2},\ldots,l_0\}\subset \mathbb{L}$ are
$\mathbb{K}$--linearly independent, since otherwise one would find a
non-zero polynomial of degree smaller than $r$ vanishing at $-d$.
Now, let  $\vec{l}_i=(l_{i,0},\ldots,l_{i,n-1})^{t}$ be the vector
of coordinates of $l_i$ in the base
$\{1,\alpha,\ldots,\alpha^{n-1}\}$. Then,
$\{\vec{1},\vec{l}_{r-2},\ldots,\vec{l}_0\}\subset {\mathbb{K}}^n$
are $\mathbb{K}$--linearly independent. Moreover, since
$(p_0(t),\ldots,p_{n-1}(t))^{t}=\vec{1}t^{r-1}+\vec{l}_{r-2}t^{r-2}+\cdots+\vec{
l}_{0}$,
there are $r$ polynomials $p_{i_j}$, $0\leq i_1<\cdots<i_r\leq n-1$,
linearly independent. By simplicity, we assume w.l.o.g. that the
first $r$ polynomials are linearly independent. Observe that this is
always possible through a permutation matrix. The new curve, that we
will continue denoting by ${\mathcal{U}}$, is not, in general, a
hypercircle. In this situation, we proceed to prove (1) and (2).

\noindent In order to prove (1), let $\mathcal{A}
\in\mathcal{M}_{n-r\times r}(\mathbb{K})$ be the matrix providing
the linear  combinations of the $n-r$ last polynomials in terms of
the first $r$ polynomials; i.e. $(p_{r}(t),\ldots,p_{n-1}(t))^t=
\mathcal{A}(p_0(t),\ldots,p_{r-1}(t))^t$. Now, given the bases
$\mathcal{B}=\{1,\ldots,t^{r-1}\}$ and $\mathcal{B}^{\star}=
\{p_0(t),\ldots,p_{r-1}(t)\}$, let $\mathcal{M}\in
\mathcal{M}_{r\times r}(\mathbb{K})$ be the transpose matrix of
change of bases from $\mathcal{B}$ to $\mathcal{B}^{\star}$.
Finally, the $n\times n$ matrix
\[
\mathcal{Q}=\left(\begin{array}{cc}\mathcal{M}&\mathcal{O}_{r,n-r}\\-\mathcal{A}
&I_{n-r}\end{array}\right)\]
defines, under the previous assumptions, the  affine transformation $\chi$.
Note that if $r=n$ then $\mathcal{Q}=\mathcal{M}$.

\noindent The proof of (2) is analogous to (1). Now, let consider
the basis $\mathcal{B}=\{1,\ldots,t^{r-1},t^r\}$ and
$\mathcal{B}^{\star}=\{p_0(t),\ldots,p_{r-1}(t), M(t)\}$. Let
$\mathcal{A} \in \mathcal{M}_{n-r\times r+1}(\mathbb{K})$ be the
matrix providing the linear combinations of the $n-r$ last
polynomials in terms of basis $\mathcal{B}^{\star}$; i.e.
$(p_{r}(t),\ldots,p_{n-1}(t))^t= \mathcal{A}
(p_0(t),\ldots,p_{r-1}(t),M(t))^t$. Let $\mathcal{M}\in
\mathcal{M}_{r+1\times r+1}(\mathbb{K})$ be the transpose matrix of
change of bases from $\mathcal{B}$ to $\mathcal{B}^{\star}$.
Finally, the $n+1\times n+1$ matrix
\[ \mathcal{Q}=\left(\begin{array}{cc}\mathcal{M}&\mathcal{O}_{r+1,n-r}\\-
\mathcal{A}&I_{n-r}\end{array}\right)\]
defines, under the previous assumptions, the  projective transformation
$\rho$.  Note that if $r=n$ then $\mathcal{Q}=\mathcal{M}$.
\end{proof}

As a direct consequence, we derive the following geometric
properties of hypercircles.

\begin{corollary}\label{corolario-propiedades-geometricas-hipercirculos-1}
In the hypothesis of Theorem \ref{teo-estructura-hipercirculo}
\begin{enumerate}
\item $\mathcal{U}$ defines a curve of degree r.
\item $\mathcal{U}$ is contained in a linear variety of dimension $r$ and it is
not
contained in a variety of dimension $r-1$.
\item $\mathcal{U}$ is a regular curve in $\mathbb{P(F)}^n$.
\item The Hilbert function of \, $\mathcal{U}$ is equal to its Hilbert
polynomial and
$h_{\mathcal{U}}(m)=mn+1$.
\end{enumerate}
\end{corollary}

\begin{proof}
All these properties are well known to hold for the {\it rational normal curve
of
degree $r$} e.g. \cite{harris}, \cite{hart}, \cite{walker}).
\end{proof}

In the following theorem, we classify the hypercircles that are affinely
equivalent over $\mathbb{K}$. We will assume that the denominator of the
generating units are not constant. The case where the units are polynomials
are described in Theorem \ref{teo-rectas}.

\begin{theorem}\label{teo-clasificacion-afin-hipercirculos}
Let ${\mathcal{U}}_i,\, i=1,2,$ be  $\alpha$-hypercircles associated to
$u_i(t)=\frac{a_it+b_i}{t+d_i}$, and let $M_i(t)$ be the minimal polynomial of
$-d_i$ over $\mathbb{K}$. Then, the following statements are equivalent:
\begin{enumerate}
\item ${\mathcal{U}}_1$ and ${\mathcal{U}}_2$ are affinely equivalent over
$\mathbb{K}$.
\item There exists a unit $\tau(t)\in\mathbb{K}(t)$ such that it maps a root
(and hence all roots) of $M_1(t)$ onto a root (resp. all roots) of $M_2(t)$.
\end{enumerate}
\end{theorem}

\begin{proof}
First of all note that, because of Theorem \ref{teo-rectas}, the
result for lines is trivial. For dealing with the general case, we
observe that, by Lemma \ref{lema-forma-reducida}, we can assume that
$u_i(t)=1/(t+d_i)$. Next, suposse that ${\mathcal{U}}_1$ and
${\mathcal{U}}_2$ are affinely equivalent over $\mathbb{K}$. By
Theorem \ref{teo-estructura-hipercirculo}, statement (1),
$[\mathbb{K}(d_1):\mathbb{K}]=[\mathbb{K}(d_2):\mathbb{K}]=r$ and
the curves ${\mathcal{U}}^{\star}_1:=\chi({\mathcal{U}}_1)$ and
${\mathcal{U}}^{\star}_2:=\chi({\mathcal{U}}_{2})$ parametrized by
$\tilde{\chi}_1(t)=(\frac{1}{M_1(t)},\ldots,\frac{t^{r-1}}{M_1(t)})$
and
$\tilde{\chi}_2(t)=(\frac{1}{M_2(t)},\ldots,\frac{t^{r-1}}{M_2(t)})$,
respectively, are affinely equivalent over $\mathbb{K}$; note that,
for simplicity we have omitted the last zero components in these
parametrizations. Therefore, there exists $\mathcal{A}=(a_{i,j})\in
GL(r,\mathbb{K})$ and $\vec{v}\in M_{r\times 1}(\mathbb{K})$, such
that $\varphi(t):= {\mathcal{A}}\,\tilde{\chi}_1(t)^t+\vec{v}$
parametrizes ${\mathcal{U}}^{\star}_2$. In consequence, since
$\varphi(t)$ and $\tilde{\chi}_2(t)$ are proper parametrizations of
the same curve, there exists a unit $\tau(t)\in\mathbb{K}(t)$ such
that $\varphi(t)=\tilde{\chi}_2(\tau(t))$. Then, considering the
first component in the above equality, one gets that
\[(a_{1,1}+\cdots+a_{1,r}t^{r-1}+v_1M_1(t))M_2(\tau(t))=M_1(t).\]
Now, substituting $t$ by $-d_1$, we obtain
\[(a_{1,1}+ \cdots+ a_{1,r} (-d_1)^{r-1} +v_1 M_1 (-d_1) ) M_2 (\tau(-d_1))
=M_1(-d_1) =0.\]
Note that $a_{1,1}+\cdots+a_{1,r}(-d_1)^{r-1}\neq 0$, because
$[\mathbb{K}(d_1):\mathbb{K}]=r$. Also, note that $\tau(-d_1)$ is well
defined, because $-d_1$ does not belong to $\mathbb{K}$. This implies that
$M_2(\tau(-d_1))=0$. So, $\tau(-d_1)$ is a root of $M_2(t)$.

Conversely, let $\tau(t)=\frac{k_1t+k_2}{k_3t+k_4}\in\mathbb{K}(t)$
be a unit that maps the root $\gamma$ of $M_1(t)$ onto the root
$\beta$ of $M_2(t)$, i.e. $\tau(\gamma)=\beta$. This relation
implies that $\mathbb{K}(\gamma)=\mathbb{K}(\beta)$ and that
$\deg{(M_1(t))}=\deg{(M_2(t))}=r$. Therefore, because of Theorem
\ref{teo-estructura-hipercirculo}, it is enough to prove that the
curves ${\mathcal{U}}^{\star}_1:=\chi({\mathcal{U}}_1)$ and
${\mathcal{U}}^{\star}_2:=\chi( {\mathcal{U}}_2)$ are  affinely
equivalent over $\mathbb{K}$.  Recall that ${\mathcal{U}}^{\star}_i$
is  parametrized by
$\varphi_i(t):=\tilde{\chi}(t)=\left(\frac{1}{M_i(t)},\ldots,\frac{t^{r-1}}{
M_i(t)}\right)$;
here again, we omit the last zero components of the parametrization.
In order to prove the result, we find an invertible matrix
${\mathcal{A}}\in GL(r,\mathbb{K})$ and a vector $\vec{v}\in
M_{r\times 1} (\mathbb{K})$, such that
${\mathcal{A}}\,\varphi^{t}_1(t)+\vec{v}= \varphi^{t}_2(\tau(t)).$
For this purpose, we consider the polynomial
$M(t)=M_2(\tau(t))(k_3t+k_4)^r\in {\mathbb{K}}[t]$. Now, since
$\tau(t)$ is a unit of $\mathbb{K}(t)$, and the roots of $M_2(t)$
are not in $\mathbb{K}$, one gets that $\deg(M)=\deg(M_2)=r$.
Moreover, since $\gamma$ is a root of $M(t)$, and taking into
account that $M_1(t)$ is the minimal polynomial of $\gamma$ over
$\mathbb{K}$ and that $\deg(M)=r=\deg(M_1)$, one has that  there
exists $c\in {\mathbb K}^*$ such that  $M(t)=cM_1(t)$. Now, in order
to determine $\mathcal{A}$ and $\vec{v}$, let us substitute
$\tau(t)$ in the $i$-th component of $\varphi_2(t)$:
\[ \frac{\tau(t)^i}{M_2(\tau(t))}=
\frac{\tau(t)^i(k_3t+k_4)^r}{M_2(\tau(t))(k_3t+k_4)^r}=\frac{
(k_1t+k_2)^i(k_3t+k_4)^{r-i}}{cM_1(t)}.\]
Since numerator and denominator in the above rational function have the
same degree, taking quotients and remainders, $\varphi_2(t)$ can be
expressed as
\[(\varphi_2(\tau(t)))_{i=1,\ldots,r}=(v_i+\frac{a_{i,1}+\cdots+a_{i,r}t^{r-1}}{
M_1(t)})_{i=1,\ldots,r},\]
for some $v_i,a_{i,j}\in {\mathbb{K}}$. Take ${\mathcal{A}}=(a_{i,j})$ and
$\vec{v}=(v_i)$.
Then, ${\mathcal{A}}(\varphi_1(t))^{t}+\vec{v}=(\varphi_2(\tau(t))^{t}$.
Finally, let us see that  ${\mathcal{A}}$ is regular. Indeed, suppose that
${\mathcal{A}}$ is
singular  and that there exists a non trivial linear relation
$\lambda_1F_1+\cdots+\lambda_{r}F_{r}=\vec{0}$, where $F_i$ denotes
the $i$-th row of $\mathcal{A}$. This implies that
$\left(\lambda_1\frac{1}{M_2(t)}+\cdots+\lambda_r\frac{t^{r-1}}{M_2(t)}
\right)\circ \tau(t)=\lambda_1v_1+\cdots+\lambda_{r}v_{r}$ is  constant,
which is impossible because
$\frac{\lambda_1+\cdots+\lambda_{r}t^{r-1}}{M_2(t)}$ is not constant and
$\tau(t)$ is a unit of $\mathbb{K}(t)$.
\end{proof}

For two true circles, there is always a real affine transformation
relating them. We have seen that this is not the case of
hypercircles. However, for algebraic extensions of degree 2 (where
the circle case fits), we recover this property for hypercircles
that are not lines.

\begin{corollary}\label{cor-clasificacion-afin-hipercirculos}
Let $\mathbb{K}(\alpha)$ be an extension of degree $2$. Then all
$\alpha$-hypercircles, that are not lines, are affinely equivalent over
$\mathbb{K}$.
\end{corollary}

\begin{proof}
By Lemma \ref{lema-forma-reducida}, we may assume that the hypercircles
are associated to units of the form $\frac{1}{t+d}$.
Now, we consider two $\alpha$-hypercircles not being lines, namely,  let
${\mathcal{U}}_i$ be the $\alpha$-hypercircle associated to $\frac{1}{t+d_i}$
for
$i=1,2$, and $d_i\not\in \mathbb{K}$. Let $d_i=\lambda_i +\mu_i \alpha$,
with $\lambda_i,\mu_i\in \mathbb{K}$ and $\mu_i\neq 0$.  Then, the unit
$\tau(t)=\tau_0+\tau_1t\in\mathbb{K}[t]$ where
$\tau_0=\frac{\mu_2\lambda_1-\mu_1\lambda_2}{\mu_1}$ and
$\tau_1=\frac{\mu_2}{\mu_1}$, verifies that $\tau(-d_1)=-d_2$. By Theorem
\ref{teo-clasificacion-afin-hipercirculos}, ${\mathcal{U}}_{1}$ and
${\mathcal{U}}_{2}$
are affinely equivalent over $\mathbb{K}$.
\end{proof}

In Corollary \ref{corolario-propiedades-geometricas-hipercirculos-1}
we have seen that the degree of a hypercircle is given by the degree
of the field extension provided by the pole of any non polynomial
generating unit. Lines are curves of degree one, a particular case
of this phenomenon. Now, we consider other kind  of hypercircles of
degree smaller than $n$. This motivates the following concept.

\begin{definition}\label{def_hipercirculo_primitivo}
Let $\mathcal{U}$ be an $\alpha$-hypercircle. If the degree of
$\mathcal{U}$ is $[\mathbb{K}(\alpha):\mathbb{K}]$, we say that it is a
primitive hypercircle. Otherwise, we say that $\mathcal{U}$ is a
non-primitive hypercircle.
\end{definition}

Regarding the  complex numbers as an extension of the reals, lines  may be
considered as circles when we define them through a Moebius transformation.
Lines are the only one curves among these such that its degree is not
$[\mathbb{C}:\mathbb{R}]$. The situation is more complicated in the
general case. Apart from lines, which have been thoroughly studied in
Theorem \ref{teo-rectas}, there are other non-primitive hypercircles. This is
not a big challenge because, as we will see, non-primitive hypercircles are
primitive on another extension. Moreover, these cases reflect some algebraic
aspects of the extension
$\mathbb{K}\subseteq\mathbb{K}(\alpha)={\mathbb{L}}$ in the geometry
of the hypercircles. Actually, we will see that there is a correspondence
between non-primitive hypercircles and the intermediate fields of
$\mathbb{K}\subseteq{\mathbb{L}}$. More precisely, let $\mathcal{U}$
be a non-primitive hypercircle associated to $u(t)=\frac{1}{t+d}$, where
$r=[\mathbb{K}(d):\mathbb{K}]<[\mathbb{L}:\mathbb{K}]=n$. In this
case, we have the algebraic extensions
$\mathbb{K}\subseteq\mathbb{K}(d)\subsetneq\mathbb{L}$. We may
consider $u(t)$ as a unit either in the extension
$\mathbb{K}\subseteq\mathbb{K}(d)$ with primitive element $d$ or in
$\mathbb{K}(d)\subsetneq\mathbb{L}$ with primitive element $\alpha$. In
the first case, $u(t)$ defines a primitive hypercircle in $\mathbb{F}^r$. In
the second case, as $u(t)$ is a $\mathbb{K}(d)$ unit, it defines a line. The
analysis of $\mathcal{U}$ can be reduced to the case of the primitive
hypercircle associated to $u(t)$ in the extension
$\mathbb{K}\subseteq\mathbb{K}(d)$.

\begin{theorem}\label{teo_paso_a_propio}
Let $\mathcal{U}$ be the non-primitive hypercircle associated to
$u(t)=\frac{at+b}{t+d}\in\mathbb{K}(\alpha)(t)$. Let $\mathcal{V}$
be the hypercircle generated by the unit $\frac{1}{t+d}$ in the
extension $\mathbb{K}\subseteq\mathbb{K}(d)$. Then, there is an
affine inclusion from $\mathbb{F}^r$ to $\mathbb{F}^n$, defined over
$\mathbb{K}$, that maps the hypercircle $\mathcal{V}$ onto
$\mathcal{U}$.
\end{theorem}

\begin{proof}
Taking into account Lemma \ref{lema-forma-reducida}, we may assume that
$u(t)=\frac{1}{t+d}$.  Let
$\phi(t)=(\phi_0(t),\ldots,\phi_{n-1}(t))\in\mathbb{K}(t)^n$ be the
parametrization of $\mathcal{U}$, obtained from $u(t)$, with respect to the
basis
$\mathcal{B}=\{1,\alpha,\ldots,\alpha^{n-1}\}$. Similarly, let
$\psi(t)=(\psi_0(t),\ldots,\psi_{r-1}(t))\in\mathbb{K}^r(t)$ be the
parametrization of the hypercircle $\mathcal{V}$, associated to $u(t)$, with
respect to the basis $\mathcal{B}^{\star}=\{1,d,\ldots,d^{r-1}\}$, where
$r=[\mathbb{K}(d):\mathbb{K}]$. The matrix $\mathcal{D}=(d_{ji})\in
\mathcal{M}_{n\times r}(\mathbb{K})$ whose columns are the coordinates of
$d^i$ with respect to $\mathcal{B}$ induces a $\mathbb{K}$-linear
transformation $\chi:\mathbb{F}^r\mapsto \mathbb{F}^n$ that maps
$\mathcal{V}$ onto $\mathcal{U}$. Indeed, as
$u(t)=\sum_{i=0}^{r-1}\psi_i(t)d^i=\sum_{j=0}^{n-1}\phi_j(t)\alpha^j$,
one has that
\[\sum_{i=0}^{r-1}\psi_i(t)d^i=\sum_{i=0}^{r-1}\psi_i(t)
\left(\sum_{j=0}^{n-1}d_{j,i}\alpha^{j}\right)=
\sum_{j=0}^{n-1}\left(\sum_{i=0}^{r-1}d_{j,i}\psi_i(t)\right)\alpha^{j}=
\sum_{j=0}^{n-1}\phi_j(t)\alpha^j.\]
Then $\phi(t)^{t}=\mathcal{D}\,\psi(t)^{t}$. Moreover, $\chi$ is one to one,
because ${\rm rank}(D)=r$.
\end{proof}

As a consequence of this theorem, every hypercircle is affinely equivalent,
over $\mathbb{K}$, to a primitive hypercircle. Therefore, the study of
hypercircles can be reduced to the study of primitives hypercircles.

\section{Properties at Infinity of a Hypercircle}
Circles have a very particular structure at infinity, namely, they
pass  through the cyclic points, i.e.  $[\pm i:1:0]$, which are
related to the minimal polynomial defining the circle as a
hypercircle as remarked in the introduction. In this section, we
will see that a similar situation occurs for more general primitive
hypercircles. More precisely, let  $\mathcal{U}$ be the primitive
hypercircle  defined by the unit $u(t)=\frac{at+b}{t+d}$. By
Corollary \ref{corolario-propiedades-geometricas-hipercirculos-1},
$\mathcal{U}$ is a parametric affine curve of degree $n$. So, there
are at most $n$ different points in the hyperplane at infinity. Let
$\phi(t)=\left(\phi_0(t),\ldots,\phi_{n-1}(t)\right)$ be the
parametrization of ${\mathcal{U}}$ generated by $u(t)$; recall that
$\phi_i(t)=\frac{p_i(t)}{M(t)}$. Thus, projective coordinates of the
points attained by $\phi(t)$ are given by
$[p_0(t):\cdots:p_{n-1}(t):M(t)]$. Now, substituting $t$ by every
conjugate $\sigma(-d)$ of $-d$, we obtain
\[[p_0(\sigma(-d)):\cdots:p_{n-1}(\sigma(-d)):0]=[
\sigma(p_0(-d)):\cdots:\sigma(p_{n-1}(-d)):0]\]
We prove next that these points are the points of the hypercircle at infinity.

\begin{lemma}\label{lema-infinito}
Let ${\mathcal{U}}$ be a primitive hypercircle associated to the unit
$u(t)=\frac{at+b}{t+d}$. The $n$ points at infinity are
\[P_j=[\sigma_j(p_0(-d)):\cdots:\sigma_j(p_{n-1}(-d)):0],\ 1\leq j\leq n\]
where $\sigma_j$ are the $\mathbb{K}$-automorphisms of the normal
closure of $\mathbb{L}=\mathbb{K}(\alpha)$ over $\mathbb{K}$.
\end{lemma}
\begin{proof}
First of all, observe that $\gcd(p_0,\ldots,p_{n-1},M)=1$, and hence
$P_j$ are well defined. Moreover, $p_i(-d)\neq 0$, for every
$i\in\{0,\ldots,n-1\}$, since $p_i(t)\in\mathbb{K}[t]$ is of degree
at most $n$ and, thus, if $p_i(-d)=0$, then
$\frac{p_i(t)}{M(t)}=c\in\mathbb{K}$ and the hypercircle would be
contained in a hyperplane. But this is impossible since
$\mathcal{U}$ is primitive (see Corollary
\ref{corolario-propiedades-geometricas-hipercirculos-1}). It remains
to prove that they are different points. Suppose that two different
tuples define the same projective point. We may suppose that
$P_1=P_j$. $P_1$ verifies that $\sum_{i=0}^{n-1}p_i(-d)\alpha^i=
(-ad+b)m(-d)\neq 0$ and $P_j$ verifies that
$\sum_{i=0}^{n-1}p_i(\sigma_j(-d))\alpha^i=(a\sigma_j(-d)+b)m(
\sigma_j(-d))=0$. Thus, $P_j$ is contained in the projective
hyperplane $\sum_{i=0}^{n-1}\alpha^iX_i=0$, but not $P_1$. Hence,
$P_1\neq P_j$.
\end{proof}

Let us check that, as in the case of circles, the points at infinity of
primitive
$\alpha$-hypercircles do not depend on the particular hypercircle.

\begin{theorem}\label{teo-infinito}
For a fixed extension $\mathbb{K}\subseteq\mathbb{K}(\alpha)$ of degree
$n$, the set of points at the infinity $P=\{P_1,\ldots,P_n\}$ of any primitive
hypercircle does not depend on the particular $\alpha$-hypercircle
$\mathcal{U}$,
but only on the algebraic extension and on the primitive element $\alpha$.
Moreover, the set $P$ is characterized by the following property:
\[\{X_0+\alpha_jX_1+\cdots+\alpha_j^{n-1}X_{n-1}=0\}\cap
\overline{\mathcal{U}}=P\setminus\{P_j\},\]
where $\alpha_j=\sigma_j(\alpha)$ are the conjugates of $\alpha$ in
$\mathbb{F}$, $1\leq j\leq n$, and $\overline{\mathcal{U}}$ is the
projective closure of $\mathcal{U}$.
\end{theorem}
\begin{proof}
Let $\mathcal{U}$ be the primitive $\alpha$-hypercircle generated be a unit
$u(t)=\frac{at+b}{t+d}$. $\overline{\mathcal{U}}$ has the projective
parametrization $[p_0(t):\cdots:p_{n-1}(t):M(t)]$. Let
$P_j=[\sigma_j(p_0(-d)):\cdots:\sigma_j(p_{n-1}(-d)):0]$. Its evaluation in the
equation of hyperplane $X_0+\alpha_kX_1+\ldots+\alpha_k^{n-1}X_{n-1}$,
yields:
\[\sum_{i=0}^{n-1}\sigma_j(p_i(-d))\alpha_k^i=
\sigma_k\left(\sum_{i=0}^{n-1}\sigma_{k}^{-1}\circ\sigma_j(p_i(-d))
\alpha^i\right)=\]
\[\sigma_k\left((a(\sigma_{k}^{-1}\circ\sigma_j(-d))+b)m(\sigma_{k}^{-1}
\circ\sigma_j(-d))\right).\]
If $j=k$, the previous expression equals
$\sigma_k\left((-ad+b)m(-d)\right)\neq 0$. If $j\neq k$, then
$\sigma_{k}^{-1}\circ\sigma_j(-d)$ is a conjugate of $-d$, different from
$-d$, because $-d$ is a primitive element. So
$m(\sigma_{k}^{-1}\circ\sigma_j(-d))=0$.

In order to show that this point does not depend on a particular hypercircle,
take the $n$ hyperplanes
$X_0+\alpha_kX_1+\cdots+\alpha_k^{n-1}X_{n-1}=0$, $k=1\ldots n$.
Every point at infinity of a hypercircle is contained in exactly $n-1$ of those
hyperplanes. Also, any of these hyperplanes contains exactly $n-1$ points at
infinity of the hypercircle. One point at infinity may be computed by solving
the linear system given by any combination of $n-1$ hyperplanes. The matrix
of the linear system is a Vandermonde matrix, each row depending on the
corresponding $\alpha_k$, so there is only one solution.
\end{proof}
\begin{remark}
Notice that this theorem provides a $n$-simplex combinatorial structure of
the points at infinity of any primitive hypercircle.
\end{remark}

The following result shows that the points at infinity can be read directly from
the minimal polynomial of $\alpha$.

\begin{proposition}\label{prop-puntos-infinito-directo}
Let $M_{\alpha}(t)$ be the minimal polynomial of $\alpha$ over $\mathbb{K}$.
Let $m_{\alpha}(t)=\frac{M_{\alpha}(t)}{t-\alpha}=
\sum_{i=0}^{n-1}l_it^i\in\mathbb{K}(\alpha)[t]$, where $l_{n-1}=1$.
Then, the points at infinity of every primitive $\alpha$-hypercircle are
$[l_0:l_1:\cdots:l_{n-2}:l_{n-1}:0]$ and its conjugates.
\end{proposition}

\begin{proof}
We consider the symmetric polynomial
$r(x,y)=\frac{M_{\alpha}(x)-M_{\alpha}(y)}{x-y}$. Substituting $(x,y)$ by
$(t,\alpha)$ we obtain that
\[r(t,\alpha)=\frac{M_{\alpha}(t)-M_{\alpha}(\alpha)}{t-\alpha}=
\frac{M_{\alpha}(t)}{t-\alpha}=m_{\alpha}(t).\] That is,
$m_{\alpha}(t)$ is symmetric in $t$ and $\alpha$. Take now the
hypercircle induced by the unit $\frac{1}{t-\alpha}=
\frac{m_{\alpha}(t)}{M_{\alpha}(t)}$. By Lemma \ref{lema-infinito},
we already know that one point at infinity is
$[p_0(\alpha):\cdots:p_{n-1}(\alpha):0],$ where
$m_{\alpha}(t)=\sum_{i=0}^{n-1}p_i(t)\alpha^i$. By symmetry,
$\sum_{i=0}^{n-1}p_i(t)\alpha^i=\sum_{i=0}^{n-1}p_i(\alpha)t^i$.
That is, $p_i(\alpha)=l_i$. Thus, the  points at infinity are
$[l_0:l_1:\cdots:l_{n-2}:1:0]$ and its conjugates.
\end{proof}

Next result deals with the tangents of a hypercircle at infinity, and it
explains
again why parabolas can not be hypercircles.

\begin{proposition}
The tangents to a primitive hypercircle at the points at infinity are not
contained in the hyperplane at infinity.
\end{proposition}

\begin{proof}
Let $\mathcal{U}$ be the primitive $\alpha$-hypercircle generated by
$\frac{at+b}{t+d}$, and $[p_0(t):\cdots:p_{n-1}(t):M(t)]$ the projective
parametrization generated by the unit. In the proof of Lemma
\ref{lema-infinito}, we have seen that $p_{n-1}(t)$ is not identically $0$,
because $p_{n-1}(-d)\neq 0$. So, we can dehomogenize w.r.t. the variable
$X_{n-1}$, obtaining the affine parametrization
$(\frac{p_0(t)}{p_{n-1}(t)},\ldots,\frac{p_{n-2}(t)}{p_{n-1}(t)},
\frac{M(t)}{p_{n-1}(t)})$ of $\mathcal{U}$ on another affine chart. We
have to check that the tangents to the curve at the intersection points with
the hyperplane $X_{n-1}=0$ are not contained in this hyperplane. The
points of $\mathcal{C}$ in the hyperplane $X_{n-1}=0$ are obtained by
substituting $t$ by $\sigma(-d)$. The last coordinate of the tangent vector is
\[\frac{M'(t)p_{n-1}(t)-M(t)p'_{n-1}(t)}{p_{n-1}(t)^2}.\]
We evaluate this expression at $\sigma(-d)$. $M(\sigma(-d))=0$ and, as all its
roots are different in $\mathbb{F}$, $M'(\sigma(-d))\neq 0$. We also know
that $\sigma(p_{n-1}(-d))\neq 0$. Hence, the last coordinate of the tangent
vector is non-zero. Thus, the tangent line is not contained in the hyperplane at
infinity.
\end{proof}

Finally, we present a property of hypercircles that can be derived from the
knowledge of its behavior at infinity. We remark a property of circles stating
that given three different points in the plane, there is exactly one circle
passing through them (which is a line if they are collinear). The result is
straightforward if we recall that there is only one conic passing throught five
points. In the case of circles, we  have the two points at infinity already
fixed,
so, given three points in the affine plane there will only be a conic (indeed a
circle if it passes through the cyclic points at infinity) through them. Even if
hypercircles are curves in n-space, surprisingly, the same occurs for
hypercircles.

We are going to prove that, given 3 different points in
$\mathbb{K}^n$, there is exactly one hypercircle passing through
them. If the points are not in general position, the resulting
hypercircle needs not to be a primitive one. First, we need a lemma
that states what are the points over $\mathbb{K}$ of the hypercircle
that are reachable by the parametrization.

\begin{lemma}\label{lemma_2_issac04}
Let $\mathcal{U}$ be the $\alpha$--hypercircle, non necessarily primitive,
associated to $u(t)=\frac{at+b}{t+d}$ with induced parametrization $\Phi(t)$.
$\Phi(\mathbb{K})=\mathcal{U}\cap\mathbb{K}^n\setminus\{\bar{a}\}$
with $a=\sum_{i=0}^{n-1}a_i\alpha^i$, $\bar{a}=(a_0,\ldots,a_{n-1}).$
\end{lemma}
\begin{proof}
We already know that $\Phi(t)$ is proper  and, obviously,
$\Phi(\mathbb{K})\subseteq\mathcal{U}\cap\mathbb{K}^n$, also,
$\bar{a}$ is not reachable by $\Phi(t)$, since otherwise one would have that
$a=u(\lambda)$ for some $\lambda$, and this implies that $ad-b=0$, which
is impossible since $u(t)$ is a unit. In order to prove the other inclusion,
write
as before $\phi_i(t)=\frac{p_i(t)}{M(t)}$, where $M(t)$ is the minimal
polynomial of $-d$ over $\mathbb{K}$. Then, we consider the ideal $I$ over
$\mathbb{F}[t,\bar{X}]$ generated by
$(p_0(t)-X_0M(t),\ldots,p_{n-1}(t)-X_{n-1}M(t))$, where
$\bar{X}=(X_0,\ldots,X_{n-1})$, and the ideal
$J=I+(ZM(t)-1)\subseteq\mathbb{F}[Z,t,\bar{X}]$. Let $I_1$ be the first
elimination ideal of $I$; i.e. $I_1=I\cap\mathbb{F}[\bar{X}]$ and let $J_2$ be
the second elimination ideal of $J$; i.e. $J_2=J\cap \mathbb{F}[\bar{X}]$.
Observe that $I\subseteq J$ and therefore $I_1\subseteq J_2$. Note that
$\mathcal{U}=V(J_2)$; i.e. $\mathcal{U}$ is the variety defined by $J_2$
over $\mathbb{F}$. Thus $\mathcal{U}\subseteq V(I_1)$. Now, let us take
$\bar{x}\in(\mathcal{U}\cap\mathbb{K}^n)\setminus \{\bar{a}\}$. Then
$\bar{x}\in V(I_1)$. Observe that, by construction, the leading coefficient of
$p_i(t)-X_iM(t)$ w.r.t. $t$ is $a_i-X_i$. Therefore, since $\bar{x}\neq
\bar{a}$ one has that at least one of the leading coefficients of the
polynomials in $I$ w.r.t. $t$ does not vanish at $\bar{x}$. Thus, applying the
Extension Theorem (see Theorem 3, pp. 117 in \cite{Cox}), there exists
$t_0\in\mathbb{F}$ such that $(t_0,\bar{x})\in V(I)$. This implies that
$p_i(t_0)-x_iM(t_0)=0$ for $i=1\ldots n-1$. Let us see that $M(t_0)\neq 0$.
Indeed, if $M(t_0)=0$ then $p_i(t_0)$ is also zero for every index and
therefore $\gcd(p_0(t),\ldots,p_{n-1}(t),M(t))\neq 1$, which is impossible.
Hence $\Phi$ is defined at $t_0$ and $\Phi(t_0)=\bar{x}$. To end up, we only
need to show that $t_0\in\mathbb{K}$. For this purpose, we note that the
inverse of $\Phi(t)$ is given by
\[P(\bar{X})=\frac{-d\sum X_i\alpha^i+b}{\sum X_i\alpha^i-a}\]
Now, since $\bar{x}\neq \bar{a}$ one deduces that $P(\bar{x})$ is well
defined, and the only parameter value generating $\bar{x}$ is
$t_0=P(\bar{x})$. Hence, the $\gcd$ of the polynomials $p_i(t)-x_iM(t)$ is a
power of $(t-t_0)$. Thus, taking into account that $p_i, M\in\mathbb{K}[t]$,
one deduces that $t_0\in\mathbb{K}$. Finally, it only remains to state that
$\bar{a}$ is generated when $t$ takes the value of the infinity of
$\mathbb{K}$. But this follows taking $\Phi(1/t)$ and substituting by $t=0$.
\end{proof}

\begin{proposition}\label{prop-3-points}
Let
$X_i=(X_{i0},\ldots,X_{i,n-1})\in\mathbb{K}^n\subseteq\mathbb{F}^n$,
$1\leq i\leq 3$ be three different points. Then, there exists only one
$\alpha$--hypercircle passing through them.
\end{proposition}
\begin{proof}
Let $Y_i=\sum_{j=0}^{n-1}X_{ij}\alpha^j\in \mathbb{K}(\alpha)$,
$1\leq i\leq 3$. Consider the following linear homogeneous system in $a,b,c,d$:
\[b=Y_1d,\ a+b=Y_2(c+d),\ a=Y_3c\]
Observe that, if the three points are different, there is only one projective
solution, namely $[a:b:c:d]$ where $a=Y_1Y_3-Y_3Y_2$, $b=Y_1Y_2-Y_1Y_3$,
$c=Y_1-Y_2$, $d=Y_2-Y_3$.

Take the unit $u(t)=\frac{at+b}{ct+d}$. It verifies that $u(0)=Y_1$,
$u(1)=Y_2$, $u(\infty)=Y_3$. Then, the hypercircle associated to $u$ passes
through $X_1, X_2, X_3$. In order to prove that this hypercircle is unique, let
$v$ be the unit associated to a hypercircle passing through the three points
and $\psi(t)$ the parametrization induced by $v(t)$. By Lemma
\ref{lemma_2_issac04}, as $X_i\in\mathbb{K}^n$, the point $X_i$ is
reached for a parameter value $t_i$ in $\mathbb{K}\cup \{\infty\}$. So,
there are three values $t_1,t_2,t_3\in\mathbb{K}\cup\{\infty\}$ such that
$v(t_i)=Y_i$. Let $\tau(t)\in\mathbb{K}(t)$ be the unique unit associated to
the transformation of the projective line $\mathbb{P(F)}$ into itself given by
$\tau(0)=t_1$, $\tau(1)=t_2$, $\tau(\infty)=t_3$. Then $v(\tau(t))=u(t)$ and
both units represents the same hypercircle.
\end{proof}

\section{Parametrization and Implicitation of a Hypercircle}

In this section, we will provide specific methods to parametrize and
implicitate hypercircles. These methods show the power of the rich
structure of hypercircles, simplifying problems that are usually
much harder in general.

Given a unit $u(t)$ defining $\mathcal{U}$, it is immediate to obtain a
parametrization of $\mathcal{U}$ (see Section 2). If $\mathcal{U}$ is given by
implicit
equations (as it is usually the case in Weil's descente method), the next
proposition shows how to parametrize it.

\begin{proposition}\label{prop-parametrizacion}
The pencil of hyperplanes
$X_0+X_1\alpha+\cdots+X_{n-1}\alpha^{n-1}=t$ parametrizes the
primitive $\alpha$--hypercircle $\mathcal{U}$.
\end{proposition}
\begin{proof}
Let $I$ be the implicit ideal of $\mathcal{U}$. Note that, since
$\mathcal{U}$ is $\mathbb{K}-$rational it is $\mathbb{K}$-definable,
and hence  a set of generators of $I$ can be taken in
$\mathbb{K}[X_0,\ldots,X_{n-1}]$. Let $u(t)$ be any unit associated
with $\mathcal{U}$ and $(\phi_0(t),\ldots,\phi_{n-1}(t))$ the
induced parametrization. Let $v(t)$ be the inverse unit of $u(t)$,
$u(v(t))=v(u(t))=t$. Then
$(\phi_0(v(t)),\ldots,\phi_{n-1}(v(t)))=(\psi_0(t),\ldots,\psi_{n-1}
(t))=\Psi(t)$
is another parametrization of $\mathcal{U}$ which is no more defined
over $\mathbb{K}$ but over $\mathbb{K}(\alpha)$. The later
parametrization is in standard form \cite{issac04}, that is
\[\sum_{i=0}^{n-1}\psi_i(t)\alpha^i=\left(\sum_{i=0}^{n-1}
\phi_i(t)\alpha^i\right)\circ v(t)=u\circ v(t)=t.\]
This implies that the pencil of hyperplanes $H_t\equiv
X_0+X_1\alpha+\cdots+X_{n-1}\alpha^{n-1}-t$ parametrizes
$\mathcal{U}$. Indeed, if $\Psi(t)$ is defined, $H_t\cap
\mathcal{U}$ consists in $n-1$ points at infinity of $\mathcal{U}$
(Theorem \ref{teo-infinito}) and $\Psi(t)$ itself. We deduce that
$\psi_i(t)-X_i$ belongs to the ideal $I+H_t$, which has a set of
generators in $\mathbb{K}(\alpha)(t)[X_0,\ldots,X_{n-1}]$. So, the
parametrization $\Psi(t)$ can be computed from $I$.
\end{proof}

Notice that the obtained parametrization $\Psi(t)$ has coefficients over
$\mathbb{K}(\alpha)$. Thus, it is not the parametrization induced by any
associated unit $u(t)$. The interest of obtaining a unit associated to a
hypercircle is that it helps us to solve the problem of reparametrizing a curve
over an optimal field extension of $\mathbb{K}$, see \cite{ARS-2}. There, it
is shown that given a parametrization $\Psi(t)\in\mathbb{K}(\alpha)^r$ of a
curve there is a hypercircle associated to it. Any unit associated to the
hypercircle reparametrizes the original curve over $\mathbb{K}$. To get a
parametrization $\phi(t)$ over $\mathbb{K}$ or, equivalently, a unit $u(t)$
associated to $\mathcal{U}$, we refer to \cite{issac04}. In addition, note that
the proof of Proposition \ref{prop-3-points} shows how to construct a unit
associated to a hypercircle, when points over ${\mathbb{K}}$ are known,
and therefore a parametrization of it.

The inverse problem, computing implicit equations of a hypercircle from the
parametrization induced by an associated unit, can be performed using classic
implicitation methods. However, the special structure of hypercircles provides
specific methods that might be more convenient.

\begin{proposition}\label{prop-inversa-por-unidad}
Let $\mathcal{U}$ be a hypercircle associated to the unit $u(t)$, and let
$v(t)$ be the inverse of $u(t)$. Let
\[v\left(\sum_{i=0}^{n-1}\alpha^i X_i\right)=
\sum_{i=0}^{n-1}\frac{r_i(X_0,\ldots,X_{n-1})}{s(X_0,\ldots,X_{n-1})}
\alpha^i,\] where $r_i,s\in {\mathbb{K}}[X_0,\ldots,X_{n-1}]$. Then,
the ideal of $\mathcal{U}$ is the elimination ideal with respect to
$Z$:
\[\mathcal{I(U)}=(r_1(\bar{X}),\ldots,r_n(\bar{X}),s(\bar{X})Z-1)\cap
\mathbb{F}[X_{0},\ldots,X_{n-1}].\]
\end{proposition}
\begin{proof}
Let  $u(t)=\frac{at+b}{t+d}$, then $v(t)=\frac{-dt+b}{t-a}$. Now, consider
\[u\left(\sum_{i=0}^{n-1}\alpha^i X_i\right)=\sum_{i=0}^{n-1}
\xi_i(X_0,\ldots,X_{n-1})\alpha^i\]
\[v\left(\sum_{i=0}^{n-1}\alpha^i X_i\right)=\sum_{i=0}^{n-1}
\eta_i(X_0,\ldots,X_{n-1})\alpha^i\]
where $\xi_i$, $\eta_j\in\mathbb{K}(X_0,\ldots,X_{n-1})$ and
$\eta_i=\frac{r_i(X_0,\ldots,X_{n-1})}{s(X_0,\ldots,X_{n-1})}$. The map
$\xi:\mathbb{F}^n\longrightarrow \mathbb{F}^n$,
$\xi=(\xi_0,\ldots,\xi_{n-1})$ is birational and its inverse is
$\eta=(\eta_0,\ldots,\eta_{n-1})$. Indeed:
\[\sum_{i=0}^{n-1}\eta_i(\xi_0(\bar{X}),\ldots,\xi_{n-1}(\bar{X}))\alpha^i=
v\left(\sum_{j=0}^{n-1}\alpha^j
\xi_j(\bar{X})\right)=\]
\[=v\left(u\left(\sum_{i=0}^{n-1}\alpha^iX_i\right)\right)=\sum_{i=0}^{n-1}
\alpha^iX_i\]
is an equality in $\mathbb{K}(\alpha)(X_0,\ldots,X_{n-1})$. We deduce that
\[\eta_i\big(\xi_0(X_0,\ldots,X_{n-1}),\ldots,\xi_{n-1}(X_0,\ldots,X_{n-1}
)\big)=X_i\]
It is clear that $\mathcal{U}$ is the image of the line $L\equiv
\{X_1=0,\ldots,X_{n-1}=0\}$ under the map $\xi$, $\mathcal{U}=\xi(L)$.
The set of points where $\xi$ is not defined is the union of the hyperplanes
$\sum_{i=0}^{n-1}\sigma_j(\alpha)^iX_i+\sigma_j(d)=0$, $1\leq j\leq n$.
The intersection of these hyperplanes with $L$ is the set of points
$(-\sigma(d)_j,0,\ldots,0)$, $1\leq j\leq n$. Thus, for a generic $p\in L$,
$\xi(p)$ is defined and belongs to $\mathcal{U}$. The result is similar for the
inverse map $\eta$. The set of points where $\eta$ is not defined is the union
of the hyperplanes $\sum_{i=0}^{n-1}\sigma_j(\alpha)^iX_i-\sigma_j(a)=0$,
$1\leq j\leq n$. These $n$ hyperplanes intersect $\mathcal{U}$ in at most
one affine point, see Proposition \ref{prop-parametrizacion}. So, for a generic
$p\in\mathcal{U}$, $\eta(p)$ is again defined and belongs to $L$. Let us
compute now the points $\bar{X}$ such that $\eta(\bar{X})$ is defined, but
it does not belong to the domain of $\xi$. If $\bar{X}$ is such a point, then
\[\sum_{i=0}^{n-1}\sigma_j(\alpha)^i\eta_i(\bar{X})+\sigma_j(d)=0.\]
As $\eta_i$ is defined over $\mathbb{K}$, applying $\sigma_j$ to the
definition of $\eta$, we obtain that
\[\sigma_j(v)\left(\sum_{i=0}^{n-1}\sigma_j(\alpha)^iX_i\right)=-\sigma_j(d)\]
But $\sigma_j(v)=\frac{-\sigma_j(d)t+\sigma_j(b)}{t-\sigma_j(a)}$. It follows
from Lemma \ref{lemma_2_issac04} that the value $-\sigma_j(d)$ cannot be
reached, even in $\mathbb{F}$. Thus, the image of $\eta$ is contained in the
domain of $\xi$.

We are ready to prove the theorem, by verifying that the set
$\mathcal{U}\setminus\{s=0\}$, which is just eliminating a finite number of
points in $\mathcal{U}$, is the set of points $\bar{X}$ such that
$r_i(\bar{X})=0$, $i\geq 1$ and $s(\bar{X})\neq 0$. If $\bar{X}\in
\mathcal{U}\setminus\{s=0\}$, then $\eta$ is defined and
$\eta(\bar{X})=(\eta_0(\bar{X}),0,\ldots,0)$. Hence
$\eta_i(\bar{X})=r_i(\bar{X})=0$. Conversely, if $\bar{X}$ is a point such
that $r_i(\bar{X})=0$ and $s(\bar{X})\neq 0$, then $\eta(\bar{X})$ is
defined and belongs to $L$. It is proven that $\xi$ is defined in
$\eta(\bar{X})$, so $\bar{X}=\xi(\eta(\bar{X}))\in\xi(L)=\mathcal{U}$. The
thesis of the theorem follows taking the Zariski closure of
$\mathcal{U}\setminus \{s=0\}$.
\end{proof}

This method  to compute the implicit equations of $\mathcal{U}$ is not free
from elimination techniques, as it  has to eliminate the variable $Z$. However,
it has the advantage that it yields already  an ideal in
$\mathbb{F}[X_{0},\ldots,X_{n-1}]$  defined over $\mathbb{K}$ and such
that it describes a non trivial variety containing the hypercircle. Namely,
$(r_1(\bar{X}),\ldots,r_{n-1}(\bar{X}))$ are polynomials over $\mathbb{K}$
whose zero set contains the hypercircle. The following example shows that the
elimination step is necessary in some cases.

\begin{example}
Let $\mathbb{Q}\subseteq\mathbb{Q}(\alpha)$ be the algebraic extension
defined by $\alpha^3+\alpha^2-3=0$. Let us consider the unit
$u(t)=\frac{(2+\alpha)t+\alpha}{t+1-\alpha}$. Its inverse is
$v(t)=\frac{(\alpha-1)t+\alpha}{t-2-\alpha}$. A parametrization of
$\mathcal{U}$ is
\[\phi(t)=\left(\frac{2t^3+6t^2+7t+3}{t^3+4t^2+5t-1},\frac{t^3+6t^2+9t+2}{
t^3+4t^2+5t-1}, \frac{t^2+4t+1}{t^3+4t^2+5t-1}\right)\]
A Gr\"obner basis of the ideal of the curve is\\
$I:=\{x_1^2-x_2x_0-x_2x_1-x_1+x_2,x_0x_1-x_2x_0-3x_2^2-2x_1+4x_2, \\
x_0^2-3x_2x_1-2x_0+2x_1+3x_2-2\}$.\\
Then, proposition \ref{prop-inversa-por-unidad} states that this ideal is
\[I=(r_1(x_0,x_1,x_2),r_2(x_0,x_1,x_2),s(x_0,x_1,x_2)Z-1)\cap\mathbb{F}[x_0
,x_1,x_2]\]
where\\
$r_1=2-8x_2+4x_2x_0+6x_2^2x_0+17x_2x_1+x_2x_0^2+3x_1-3x_1^2x_2+x_0^3-x_0^2x_1+4x
_0x_1-12x_2^2-8x_1^2+9x_2^3+3x_1^3-3x_0^2-9x_0x_1x_2,$\\
$r_2=-2-7x_2+4x_2x_0-x_2x_1+8x_1-2x_0-2x_0x_1+6x_2^2-2x_1^2+x_0^2,$\\
$s=9x_2^3+6x_2^2x_0-12x_2^2+5x_2x_0-17x_2-3x_1^2x_2-9x_0x_1x_2+x_2x_0^2+24x_2x_1
+3x_1^3+8x_0+4x_0x_1-5x_0^2-x_0^2x_1+5x_1-9x_1^2-7+x_0^3$.\\
But, if we take $J=(r_1,r_2)$, then $J\subsetneq I$. The saturation of $J$ with
respect to $I$ is $J:I^\infty=(x_1^2-x_0x_2-x_1x_2-2x_1+3x_2+1,
x_0x_1-x_0x_2-3x_2^2-x_0-2x_1+2x_2+2, x_0^2-3x_1x_2-4x_0+3x_2+4)$\\
This ideal corresponds to the union of the line
\[\left\{\begin{matrix}-\alpha x_0&&+3x_2&=&-2\alpha\\
(\alpha+\alpha^2)x_0&-3x_1&&=&-3+2\alpha+2\alpha^2\end{matrix}\right.\]
and its conjugates.
\end{example}

Next theorem shows an alternative method to implicitate a hypercircle without
using any elimination techniques. It is based on properties of the normal
rational curve of degree $n$.

\begin{theorem}\label{teo-inversa-por-la-normal}
Let
$\varphi(t)=(\frac{q_0(t)}{N(t)},\ldots,\frac{q_{n-1}(t)}{N(t)})$ be
a proper parametrization of a primitive hypercircle $\mathcal{U}$
with coefficients in $\mathbb{F}$. Let $I$ be the homogeneous ideal
of the rational normal curve of degree $n$ in $\mathbb{P(F)}^n$
given by a set of homogeneous generators
$h_1(\bar{Y}),\ldots,h_r(\bar{Y})$. Let $\mathcal{Q}\in
\mathcal{M}_{n+1\times n+1}(\mathbb{F})$ be the matrix that carries
$\{q_0(t),\ldots,q_{n-1}(t),N(t)\}$ onto $\{1,t,\ldots,t^n\}$. Let
\[f_i(\bar{X})=h_i\left(\sum_{j=0}^{n}\mathcal{Q}_{0j}X_j,\ldots,\sum_{j=0}^{n}
\mathcal{Q}_{nj}X_j\right),\ 1\leq i\leq r.\]
Then $\{f_1,\ldots,f_r\}$ is a set of generators of the homogeneous ideal of
$\mathcal{U}$.
\end{theorem}

\begin{proof}
If the parametrization is proper,
$\{q_0(t),\ldots,q_{n-1}(t),N(t)\}$ is a basis of the polynomials of
degree at most $n$. This follows from the fact shown in Corollary
\ref{corolario-propiedades-geometricas-hipercirculos-1} that a
primitive hypercircle is not contained in any hyperplane. Note that
a projective point $\bar{X}$ belongs to $\mathcal{U}$ if and only if
$\mathcal{Q}(\bar{X})$ belongs to the rational normal curve, if and
only if $h_i(\mathcal{Q}(\bar{X}))=0$, $1\leq i\leq r$.
\end{proof}

\begin{remark}\mbox{}
\begin{itemize}
\item It is well known that the set of polynomials $\{Y_iY_{j-1}-Y_{i-1}Y_j\ |\
1\leq i,j\leq n\}$ is a generator set of $I$ (see \cite{harris}).
\item Notice that it is straightforward to compute $Q$ from the
parametrization. Therefore, we have an effective method to compute the
implicit ideal of the projective closure of $\mathcal{U}$. The affine ideal of
$\mathcal{U}$ can be obtained by dehomogenization $X_n=1$.
\item If the parametrization is given by polynomials over an algebraic
extension $\mathbb{K}(\beta)$ of $\mathbb{K}$, then the coefficients of
$f_i$ belongs to $\mathbb{K}(\beta)$. Moreover, if we write
$f_i(\bar{X})=\sum_{j=0}^{m}f_{ij}(\bar{X})\beta^j$, with
$f_{ij}\in\mathbb{K}[\bar{X}]$, then, $\{f_{ij}\}$ is a set of generators over
$\mathbb{K}$ of the hypercircle $\mathcal{U}$.
\item In practice, this method is much more suited to compute an implicitation
of a hypercircle than the method presented in Proposition
\ref{prop-inversa-por-unidad}.
\end{itemize}
\end{remark}

\begin{example}
The implicit equations of a hypercircle can be computed by classical
implizitation methods, for example Gr\"obner basis or with the two methods
presented in Proposition \ref{prop-inversa-por-unidad} and Theorem
\ref{teo-inversa-por-la-normal}. Here, we present two cases that show the
practical behavior of these methods. The first example considers the algebraic
extension $\mathbb{Q}\subseteq\mathbb{Q}(\alpha)$, where
$\alpha^4+\alpha^2-3$ and the unit
$u=\frac{(1-\alpha^3)t+\alpha^2}{t+1+2\alpha-3\alpha^2}$. The
parametrization of the hypercircle is given by
\[\phi_0=\frac{t^4+15t^3+22t^2+101t-195}{t^4+10t^3-17t^2-366t+233},
\phi_1=\frac{-11t^3-73t^2+65t-114}{t^4+10t^3-17t^2-366t+233},\]
\[\phi_2=\frac{2t^3+57t^2-25t-59}{t^4+10t^3-17t^2-366t+233},
\phi_3=\frac{-t^4-6t^3+4t^2+17t-56}{t^4+10t^3-17t^2-366t+233}.\]
The second example starts from the extension
$\mathbb{Q}\subseteq\mathbb{Q}(\beta)$, where $\beta$ is such that
$\beta^4+3\beta+1=0$. Here, the unit defining $\mathcal{U}$ is
$u=\frac{(1+\beta-\beta^2)t+1+\beta^3}{t+1+\beta^2-\beta^3}$ and the
parametrization induced by $u(t)$ is
\[\psi_0=\frac{t^4+11t^3+47t^2+95t+72}{t^4+13t^3+62t^2+126t+81},
\psi_1=\frac{t^4+7t^3+15t^2+17t+9}{t^4+13t^3+62t^2+126t+81},\]
\[\psi_2=\frac{-t^4-10t^3-31t^2-23t}{t^4+13t^3+62t^2+126t+81},
\psi_3=\frac{t^3+13t^2+42t+36}{t^4+13t^3+62t^2+126t+81}.\] The
running times for computing the implicit ideal (using a Mac Xserver
with 2 processors G5 2.3 GHz, 2 Gb RAM Maple 10) are
 \begin{center}\begin{tabular}{|l|r@{.}l|r@{.}l|}\hline
 &\multicolumn{2}{c|}{Example 1}&\multicolumn{2}{c|}{Example 2}\\\hline
Gr\"obner basis method&0&411&0&332\\
Proposition \ref{prop-inversa-por-unidad}&2&094&2&142\\
Theorem \ref{teo-inversa-por-la-normal}&0&059&0&021\\\hline
 \end{tabular}\end{center}

We refer the interested reader to \cite{RSTV} for a brief discussion and
comparison of the running times of these algorithms.
\end{example}

\section{Characterization of Hypercircles}
In the introduction, we defined algebraically a circle as the conic
such that its homogeneous part is $x^2+y^2$ and contains an infinite
number of real points. The condition on the homogeneous part is
equivalent to impose that the curve passes through the points at
infinity $[\pm i:1:0]$. Analogously, hypercircles are regular curves
of degree $n$ with infinite points over the base field passing
through the points at infinity described in Theorem
\ref{teo-infinito}. The following result shows that this is a
characterization of these curves.
\begin{theorem}\label{teorema-caracterizacion} Let
$\mathcal{U}\subseteq\mathbb{F}^n$ be an algebraic set of degree $n$
such that all whose components are of dimension 1. Then, it is a
primitive $\alpha$-hypercircle if and only if it has an infinite
number of points with coordinates in $\mathbb{K}$ and passes through
the set of points at infinity characterized  in Theorem
\ref{teo-infinito}.
\end{theorem}
\begin{proof}
The only if implication is trivial. For the other one, let
$\mathcal{U}\subseteq\mathbb{F}^n$ be an algebraic set of pure
dimension 1 and degree $n$ passing through $P=\{P_1,\ldots,P_n\}$,
the $n$ points at infinity of a primitive $\alpha$-hypercircle.
Suposse that $\mathcal{U}$ has infinite points with coordinates in
$\mathbb{K}$. Then, we are going to prove that $\mathcal{U}$ is
irreducible. Let $\mathcal{W}$ be an irreducible component of
$\mathcal{U}$ with infinite points in $\mathbb{K}$. Note that, since
$\mathcal{W}$ is irreducible and contains infinitely many points
over $\mathbb{K}$, the ideal $\mathcal{I(W)}$ over $\mathbb{F}$ is
generated by polynomials over $\mathbb{K}$ (see Lemma  2 in
\cite{ARS-1}). Let $q$ be any point at infinity of $\mathcal{W}$;
then $q\in P$. As $\mathcal{W}$ is $\mathbb{K}$-definable it follows
that $\mathcal{W}$ also contains all conjugates of $q$. Thus, $P$ is
contained in the set of points at infinity of $\mathcal{W}$. It
follows that $\mathcal{W}$ is of degree at least $n$; since
$\mathcal{W}\subseteq\mathcal{U}$, $\mathcal{U}=\mathcal{W}$.
Therefore, $\mathcal{U}$ is irreducible and $\mathcal{I(U)}$ is
generated by polynomials with coefficients over $\mathbb{K}$. Now,
consider the pencil of hyperplanes $H_t\equiv
X_0+X_1\alpha+\cdots+X_{n-1}\alpha^{n-1}-t$, where $t$ takes values
in $\mathbb{F}$. Notice that $\overline{H_t}\cap
P=\{P_2,\ldots,P_n\}$. Thus, $P_1\in
\overline{\mathcal{U}}\setminus\overline{H_t}$ so, for all $t$,
$\mathcal{U}\not\subseteq H_t$. Moreover, for every point
$p=(p_0,\ldots,p_{n-1})\in \mathcal{U}$,
$t(p)=\sum_{i=0}^{n-1}p_i\alpha^i\in\mathbb{F}$ is such that
$\overline{H}_{t(p)}\cap
\overline{\mathcal{U}}=\{p,P_2,\ldots,P_{n}\}$. The cardinal of
$\{t(p)\ |\ t\in \mathcal{U}\}$ is infinite, since otherwise, by the
irreducibility of $\mathcal{U}$,  it would imply that there is a
$t_0$ such that $\mathcal{U}\subseteq H_{t_0}$, which is impossible.
So, for generic $t$, the intersection is
$\overline{H}_t\cap\overline{\mathcal{U}}=\{p(t),P_2,\ldots,P_{n}\}.$
Let us check that the coordinates of $p(t)$ are rational functions
in $\mathbb{K}(\alpha)(t)$. Take the ideal $\mathcal{I(U)}$ of
$\mathcal{U}$. The ideal of $p(t)$ (as a point in $\mathbb{F}(t)^n$)
is $I+H_t$, defined over $\mathbb{K}(\alpha)(t)$. The reduced
Gr\"obner basis of the radical $I+H_t$ is of this kind
$(X_0-\psi_0,\ldots,X_{n-1}-\psi_{n-1})$ and it is also defined over
$\mathbb{K}(\alpha)(t)[X_0,\ldots,X_{n-1}]$. Hence,
$(\psi_0,\ldots,\psi_{n-1})$ is a
$\mathbb{K}(\alpha)$-parametrization of $\mathcal{U}.$ Thus, since
$\mathcal{U}$ is irreducible, it is rational. Moreover
$\sum_{i=0}^{n-1}(\psi_i(t))\alpha^i=t$ and the parametrization is
proper. As the curve is rational and has an infinite number of
points over $\mathbb{K}$, it is parametrizable over $\mathbb{K}$ (it
follows, for example from the results in \cite{SW-Symbolic}). Let
$u(t)$ be a unit such that $\Psi\circ
u(t)=(\phi_0(t),\ldots,\phi_{n-1}(t))$ is a parametrization over
$\mathbb{K}$, where $\phi_i(t)\in\mathbb{K}(t)$ and
$\sum_{i=0}^{n-1}\phi_i(t)\alpha^i=u(t)$. We conclude that
$\mathcal{U}$ is the hypercircle associated to the unit $u(t)$.
\end{proof}

Remark that a parametric curve, definable over $\mathbb{K}$ and with a regular
point
over $\mathbb{K}$, is parametrizable over the same field; for this, it is enough
to $\mathbb{K}$-birationally project the curve over a plane, such that the
$\mathbb{K}$-regular point stays regular on the projection, and then apply the
results in \cite{SW-Symbolic}. Then, a small modification of the proof above,
yields the
following:

\begin{theorem}\label{teorema-caracterizacion-dos} Let
$\mathcal{U}\subseteq\mathbb{F}^n$ be a 1-dimensional irreducible algebraic set
of degree $n$, definable over $\mathbb{K}$ .
Then, it is a
primitive $\alpha$-hypercircle if and only if it has a regular point with
coordinates in $\mathbb{K}$ and passes through
the set of points at infinity characterized  in Theorem
\ref{teo-infinito}.
\end{theorem}

\section{An Application}
As mentioned in the introduction, hypercircles play an important
role in the problem of the optimal-algebraic reparametrization of a
rational curve (see \cite{ARS-1}, \cite{ARS-2}, \cite{issac04}
\cite{SV-Polynomial}, \cite{SV-Quasipolynomial} for further
details). Roughly speaking, the problem is as follows. Given a
rational $\mathbb{K}$--definable curve $\mathcal C$ by means of a
proper rational parametrization over $\mathbb{K}(\alpha)$, decide
whether $\mathcal C$ can be parametrized over $\mathbb{K}$ and, in
the affirmative case, find a change of parameter transforming the
original parametrization into a parametrization over $\mathbb{K}$.
In \cite{ARS-2}, a $\mathbb{K}$--definable algebraic variety in
${\mathbb{F}}^{n}$, where $n=[{\mathbb{K}}(\alpha): {\mathbb{K}}]$,
is associated to $\mathcal C$. This variety is called the associated
Weil (descente) parametric variety. In \cite{ARS-2}, it is proved
that this Weil variety has exactly one one-dimensional component iff
$\mathcal{C}$ is $\mathbb{K}$--definable (which is our case) and, in
this case, $\mathcal{C}$ can be parametrized over $\mathbb{K}$ iff
this one-dimensional component is a hypercircle. Moreover, if it is
a hypercircle a proper rational parametrization over $\mathbb{K}$ of
the hypercircle generates the change of parameter one is looking
for; namely its generating unit.

In the following example, we illustrate how to use the knowledge of the
geometry of hypercircles to help solving the problem.
Suppose given the parametric curve
\[\mathcal{C}\simeq (\eta_1(t),\eta_2(t))=\]
\[\left(
\frac{(-2t^4-2t^3)\alpha-2t^4}{6\alpha^2t^2+(4t^3-2)\alpha+t^4-8t},
\frac{ -2t^4\alpha}{6\alpha^2t^2+(4t^3-2)\alpha+t^4-8t}\right)\]
where $\alpha$ is algebraic over $\mathbb{Q}$ with minimal polynomial $x^3+2$.
We follow  Weil's descente method presented in \cite{ARS-2} to associate a
hypercircle to $\mathcal{C}$.
The method consists in writing $\eta_i(\sum_{j=0}^2
t_j\alpha^j)=\sum_{j=0}^2 \frac{q_{ij}(t_0,t_1,t_2)}{N(t_0,t_1,t_2)}$.
In this situation $\mathcal{C}$ is $\mathbb{Q}-$definable if and only if
\[\mathcal{V}=\overline{V(q_{11},q_{12},q_{21},q_{22})\setminus V(N)}\]
is of dimension 1.
Moreover, $\mathcal{C}$ is $\mathbb{Q}$-parametrizable if and only if the
one-dimensional component of $\mathcal{V}$ is an  $\alpha$-hypercircle.
For this example, the equations of $\mathcal{V}$ are:\\
$\mathcal{V}=V(
2t_0^3t_2-4t_2^4+3t_0^2t_1^2+2t_1^3t_2+2t_0t_2^2+2t_1^2t_2-t_0^2t_1+6t_0t_1t_2^2
,
-6t_0^2t_1t_2+t_0^4+2t_0t_1^2-8t_0t_2^3-2t_0t_1^3+2t_0^2t_2-4t_1t_2^2-12t_1^2t_2
^2,
12t_2^2t_1^3-9t_0t_1t_2^3+6t_2^5-4t_0t_1^3-2t_0^2t_1t_2+4t_1^2t_2^2-4t_0t_2^3,
9t_0t_1^2t_2^2-9t_0^2t_2^3-2t_0^3t_2-2t_1^3t_2+6t_0t_1t_2^2-2t_2^4+t_0^2t_1-2t_1
^2t_2-2t_0t_2^2,
6t_0^2t_1t_2^2+12t_1^2t_2^3-t_0^3t_1-2t_0t_1^2t_2-2t_0^2t_2^2+8t_1t_2^3,
6t_0^3t_2^2+9t_0t_1t_2^3-6t_2^5+2t_0t_1^3-2t_0^2t_1t_2+4t_1^2t_2^2+8t_0t_2^3,
18t_2t_1^4+36t_2^4t_1+14t_0^3t_2+32t_1^3t_2+12t_0t_1t_2^2-4t_2^4-7t_0^2t_1+14t_1
^2t_2+14t_0t_2^2,
6t_0t_1^3t_2+2t_0t_1^2t_2+t_0^3t_1+2t_0^2t_2^2-8t_1t_2^3+12t_2^4t_0,
9t_0^3t_2t_1-36t_2^4t_1-4t_0^3t_2-4t_1^3t_2+12t_0t_1t_2^2-4t_2^4+2t_0^2t_1-4t_1^
2t_2-4t_0t_2^2,
6t_1^5+48t_1^2t_2^3-36t_2^4t_0-11t_0^3t_1+6t_1^4+14t_0t_1^2t_2-22t_0^2t_2^2+64t_
1t_2^3,
3t_1^4t_0+6t_0t_1t_2^3+2t_0t_1^3+t_0^2t_1t_2-2t_1^2t_2^2+2t_0t_2^3,
27t_2^4t_1^2-27t_0t_2^5-9t_0^2t_2^3+9t_2^4t_1-2t_0^3t_2-2t_1^3t_2+6t_0t_1t_2^2-2
t_2^4+t_0^2t_1-2t_1^2t_2-2t_0t_2^2, 6t_2^4t_0^2+12t_2^5t_1-5t_0t_1t_2^3+2t_2^5,
t_0t_2^5t_1+2t_2^7)$

Thus the main point is to verify that this curve is a hypercircle. If
$\mathcal{V}$ is
a hypercircle, then its points at infinity must be as in Theorem 6.1. So, let us
first of
all check whether this is the case.
The set of generators of the defining ideal  form a Gr\"obner basis with respect
 to a
graded order, thus  to compute the points at infinity we take the set of
leading forms of these polynomials.

Leading forms$=\{t_0^4-2t_0t_1^3-6t_0^2t_1t_2-12t_1^2t_2^2-8t_0t_2^3,
2t_0^3t_2-4t_2^4+3t_0^2t_1^2+2t_1^3t_2+6t_0t_1t_2^2,
9t_0t_1^2t_2^2-9t_0^2t_2^3,
12t_2^2t_1^3-9t_0t_1t_2^3+6t_2^5,
6t_0^2t_1t_2^2+12t_1^2t_2^3,
6t_0^3t_2^2+9t_0t_1t_2^3-6t_2^5,
18t_2t_1^4+36t_2^4t_1,
t_0t_2^5t_1+2t_2^7,
6t_0t_1^3t_2+12t_2^4t_0,
9t_0^3t_2t_1-36t_2^4t_1,
6t_1^5+48t_1^2t_2^3-36t_2^4t_0,
3t_1^4t_0+6t_0t_1t_2^3,
27t_2^4t_1^2-27t_0t_2^5,
6t_2^4t_0^2+12t_2^5t_1\}$

The solutions of this system, after dehomogenizing $\{t_2=1\}$, are
$t_0=t_1^2,t_1^3+2=0$. That is, the points at infinity are of the
form $[\alpha^2_i:\alpha_i:1:0]$, $\frac{x^3+2}{x-\alpha}=x^2+\alpha
x+\alpha^2$. Thus, by Proposition
\ref{prop-puntos-infinito-directo}, the points at infinity of
$\mathcal{V}$ remind those of an $\alpha$-hypercircle.

Now, following Proposition \ref{prop-parametrizacion}, we may try to parametrize
$\mathcal{V}$ by the pencil of hyperplanes $t_0+\alpha t_1+\alpha^2 t_2-t$.
Doing so, we obtain the parametrization
\[\left(\frac{(\alpha^2+2\alpha t+t^2)t}{3\alpha t+\alpha^2+3t^2},
\frac{ -1/2\alpha^2t^3}{3\alpha t+\alpha^2+3t^2}, \frac{-1/2\alpha
t^2(t+\alpha)}{3\alpha t+\alpha^2+3t^2}\right).\]
Remark that this parametrization can also be computed by means of inverse
computation techniques as described in \cite{SV-Quasipolynomial}. Then, by
direct
computation, we observe that the parametric irreducible curve defined by this
parametrization is of degree 3, passes through the point $(0,0,0)$ and this
point is
regular. Moreover, it is $\mathbb{Q}$-definable, since it is the only
1-dimensional
component of $\mathcal{V}$ (see \cite{ARS-2}), which is, by construction, a
$\mathbb{Q}$-definable variety. It follows from
\ref{teorema-caracterizacion-dos} that
it is a hypercircle.

Then, from this parametrization, the algorithm presented in \cite{issac04}
computes a
unit $u(t)=\frac{2}{2t+\alpha^2}$ associated to $\mathcal{V}$. So, $\mathcal{V}$
is
the hypercircle associated to $u(t)$ and $\mathcal{C}$ is parametrizable over
$\mathbb{Q}$. In particular, the parametrization of $\mathcal{V}$ associated to
$u(t)$
is $\left(\frac{2t^2}{2t^3+1}, \frac{-1}{2t^3+1}, \frac{-t}{2t^3+1}\right)$.
Moreover,
the unit $u(t)$ gives the change of parameter we need to compute a
parametrization of
$\mathcal{C}$ over the base field (see \cite{ARS-2}), namely:
\[\eta\left(u(t)\right)=\left(\frac{t+1}{t^4}, \frac{1}{t^4}\right).\]

\end{document}